\documentclass[11pt,reqno]{amsart}
\usepackage{verbatim}
\usepackage{amssymb}
\usepackage{mathrsfs,amsfonts,amsmath,amssymb,epsfig,amscd,xy,amsthm,pb-diagram} 
\usepackage{tikz-cd}
\usepackage{geometry}
\usepackage[pdftex,hyperindex]{hyperref}
\usepackage{tikz}
\usetikzlibrary{math}

\numberwithin{equation}{section}       

\newtheorem{theorem}{Theorem}[section]
\newtheorem{proposition}[theorem]{Proposition}

\newtheorem{lemma}[theorem]{Lemma}

\newtheorem{corollary}[theorem]{Corollary}
\newtheorem{question}[theorem]{Question}
\newtheorem{definition}[theorem]{Definition}

\newtheorem*{mainthm}{Main Theorem}
\newtheorem*{ackn}{Acknowledgment}

\theoremstyle{plain}

\theoremstyle{remark}

\newtheorem{remark}[theorem]{Remark}

\newcommand{\C}{{\mathbb C}}

\renewcommand{\P}{{\mathbb P}}
\newcommand{\Q}{{\mathbb Q}}

\newcommand{\R}{{\mathbb R}}
\newcommand{\Z}{{\mathbb Z}}
\newcommand{\N}{{\mathbb N}}

\newcommand{\cC}{{\mathcal C}}

\DeclareMathOperator{\mymod}{mod}

\DeclareMathOperator{\im}{Im}
\DeclareMathOperator{\re}{Re}

\DeclareMathOperator{\Pic}{Pic}

\DeclareMathOperator{\ord}{ord}

\DeclareMathOperator{\trop}{trop}

\DeclareMathOperator{\bN}{\mathbb{N}}

\DeclareMathOperator{\NS}{NS}
\renewcommand{\Re}{{\operatorname{Re}}}
\renewcommand{\Im}{{\operatorname{Im}}}

\newcommand{\bP}{{\mathbb P}}
\newcommand{\bZ}{{\mathbb Z}}

\newcommand{\bfx}{{\mathbf{x}}}

\newcommand{\bC}{{\mathbb C}}
\newcommand{\bR}{{\mathbb R}}
\newcommand{\bA}{{\mathbb A}}
\newcommand{\bQ}{{\mathbb Q}}

\newcommand{\cO}{\mathcal{O}}

\newcommand{\cP}{\mathcal{P}}

\DeclareMathOperator{\charac}{char}
\newcommand{\bk}{\Bbbk}

\newcommand{\cCtor}{\cC^{\mathrm{tor}}}
\newcommand{\cPtor}{\cP^{\mathrm{tor}}}
\newcommand{\pitor}{\pi_{\mathrm{tor}}}
\newcommand{\imunit}{\mathsf{i}}

\newcommand{\Linmap}{\Lambda}
\newcommand{\eg}{{\rm e.g.\ }} 
\newcommand{\ie}{{\rm i.e.\ }}

\newcommand{\tto}{\dashrightarrow}

\title{A transcendental dynamical degree}
\date{\today}

\author{Jason P. Bell \and Jeffrey Diller \and Mattias Jonsson}

\address{Department of Pure Mathematics\\
  University of Waterloo\\
  Waterloo, ON\\
  Canada N2L 3G1}
\email{jpbell@uwaterloo.ca}

\address{Department of Mathematics\\
  University of Notre Dame\\
  Notre Dame, IN 46556\\
  USA}
\email{diller.1@nd.edu}

\address{Dept of Mathematics\\
  University of Michigan\\
  Ann Arbor, MI 48109-1043\\
  USA}
\email{mattiasj@umich.edu}


\begin{document}
\begin{abstract}
  We give an example of a dominant rational selfmap of the projective plane whose
  dynamical degree is a transcendental number.
\end{abstract}
	
\maketitle

%
%
%
%
\section*{Introduction}
The most fundamental dynamical invariant of a dominant rational selfmap $f\colon X\dashrightarrow X$
of a smooth projective variety is, arguably, its (first) \emph{dynamical degree} $\lambda(f)$. It can be defined, using intersection numbers, as $\lim_{n\to\infty}(f^{n*}H\cdot H^{\dim X-1})^{1/n}$, where $H$ is any ample divisor. The limit does not depend on the choice of $H$, and it is invariant under birational conjugacy: if $h\colon X'\dashrightarrow X$ is a birational map, then $f':=h^{-1}\circ f\circ h\colon X'\dashrightarrow X'$ is a dominant rational map with $\lambda(f')=\lambda(f)$.  

The dynamical degree is often difficult to compute.  If $f$ is \emph{algebraically stable} in the sense that $f^{n*}=f^{*n}$ for the induced pullbacks of divisors on $X$~\cite{Sib99}, then $\lambda(f)$ is equal to the spectral radius of the $\bZ$-linear operator $f^*\colon\NS_\bR(X)\to\NS_\bR(X)$ on the real N\'eron-Severi group $\NS_\bR(X) := \NS(X)\otimes_\bZ \bR$; hence $\lambda(f)$ is an algebraic integer in that case.  For certain classes of maps, such as birational maps of $\bP^2$~\cite{DF01} or polynomial maps of $\bA^2$~\cite{FJ07,FJ11}, we can achieve algebraic stability after birational conjugation; hence the dynamical degree is an algebraic integer in these cases.  It has been shown, moreover, that the set of dynamical degrees of all rational maps (algebraically stable or not, and over all fields) is countable~\cite{BF00,Ure18}.  

All of this leads naturally to the question~\cite[Conjecture 13.17]{BIJMST19}: is the dynamical degree \emph{always} an algebraic integer, or at least an algebraic number? 
Surprisingly,  the answer is negative:

\begin{mainthm}
Let $\bk$ be a field with $\mathrm{char}(\bk) \neq 2$.  Then there exists a dominant rational map $f\colon\bP^2_\bk\dashrightarrow\bP^2_\bk$ whose dynamical degree is a transcendental number.
\end{mainthm}

Our examples are completely explicit, of the form $f=g\circ h$, where 
\begin{equation*}
g(y_1,y_2) = \left(-y_1\frac{1-y_1+y_2}{1-y_1-y_2},-y_2\frac{1+y_1-y_2}{1-y_1-y_2}\right)
\end{equation*}
is a fixed birational involution, conjugate by a projective linear map to the standard Cremona involution $(y_1,y_2) \mapsto (y_1^{-1},y_2^{-1})$, and $h(y_1,y_2) = (y_1^ay_2^b,y_1^{-b}y_2^a)$ is a monomial map. We show that if $(a+b\imunit)^n \notin \R$ for all integers $n>0$, then $\lambda(f)$ is transcendental.  Favre~\cite{Fav03} showed that under the same condition on $a+b\imunit$, the monomial map $h$ cannot be birationally conjugated to an algebraically stable map, though $\lambda(h) = |a+b\imunit|$ is still just a quadratic integer.  Rational surface maps, such as $f$, that preserve a rational 2-form were considered as a class by the second author and J.-L. Lin in~\cite{DL16} (see also~\cite{Bla13}) where it was shown that failure of stabilizability for $h$ implies the same for $f$.  Note that the restriction $\mathrm{char}(\bk)\neq 2$ is needed only to ensure that $g$ is non-trivial.  

\subsection*{Strategy of the proof}
Our first step toward showing $\lambda(f)$ is transcendental is to relate degrees of iterates of $h$ to those of $f$.  Writing $d_j:=\deg h^j:=(h^{j*}H\cdot H)$ for $j\ge 0$, with $H\subset\bP^2$ a line, we show in \S 1-2 that the dynamical degree $\lambda=\lambda(f)$ is the unique positive solution to the equation
\begin{equation*}
  \sum_{j=1}^\infty d_j\lambda^{-j}=1.\tag{$\star$}
\end{equation*}
In order to derive~($\star$), it is useful to consider the lift of $f$ to various blowups of $\bP^2$.  We use the language of $b$-divisors to coordinate information about divisors in different blowups.  These transform naturally and functorially under the maps $h$ and $g$, so the additional terminology is convenient for understanding the degree growth of $f$, see~\cite{FJ07,BFJ08,Can11,FJ11}.  Here we make use of the additional fact that $h$ and $g$ interact well with the toric structure of $\bP^2$. This is of course clear for the monomial map $h$, but less evident for the involution $g$.

One computes by elementary means that $d_j = \re (\gamma(j)\zeta^j)$, where $\zeta = a+b\imunit$ and $\gamma(j)\in\{-2,\pm 2\imunit, 1\pm 2\imunit\}$ is chosen to be whichever element maximizes the right side.  The condition $\zeta^n \notin \R$ means that the argument of $\zeta$ is
$$
\operatorname{Arg}(a+b\imunit) = 2\pi \theta,
$$ 
for $\theta\in (0,1)$ irrational.  Were $\theta$ rational, the Gaussian integer $\gamma(j)$ would be periodic in $j$, the analytic function
\begin{equation*}
  \Delta_h(z) := \sum_{j\geq 1} d_j z^j
\end{equation*}
rational, and $\lambda$ algebraic.  However, as Hasselblatt--Propp~\cite{HP07} observed, when $\theta$ is irrational, the sequence $(d_j)_{j\ge1}$ does not satisfy any linear recurrence relation.  

One therefore suspects that $\Delta_h(\alpha)$ is unlikely to be algebraic for any given algebraic number $\alpha\ne0$ in the domain of convergence for the series; in particular $\Delta_h(1/\lambda) = 1$ should force $\lambda$ to be transcendental.  There are many results of this type in the literature, see~\eg~\cite{Nis96,FM97,AC03,AC06,Beu06,AB07a,AB07b,BBC15}, but we were not able to locate one that implies directly that at least one of $\lambda$ and $\Delta_h(\lambda^{-1})$ must be transcendental.  Instead, we present in \S 3 a proof based on results by Evertse and others on $S$-unit equations, see~\cite{EG}; these in turn rely on the $p$-adic Subspace Theorem by Schlickewei~\cite{Sch77}. Our method draws inspiration from earlier work of Corvaja and Zannier~\cite{CZ02} and Adamczewski and Bugeaud~\cite{AB07a,AB07b}, who used the subspace theorem to establish transcendence of special values of certain classes of power series.  

The idea is that if $m/n$ is a continued fraction approximant of $\theta$, then
$\zeta^n$ is nearly real, the Gaussian integers $\gamma(j)$ are nearly $n$-periodic in $j$, and $\Delta_h(z)$ is well-approximated by the rational function $\Delta_h^{(n)}(z)=(1-z^n)^{-1}\sum_{j=1}^nd_jz^j$ obtained by assuming the $\gamma(j)$ are precisely $n$-periodic.  If the approximations improve sufficiently quickly with $n$ and $\alpha$ is algebraic, then $\Delta^{(n)}_h(\alpha)$ approximates $\Delta_h(\alpha)$ too well for the latter to also be algebraic.  Unfortunately this seems a little too much to hope for without knowing more about how well $\theta$ agrees with its approximants.

To deal with the possibility that $\theta$ is badly approximable by rational numbers, we need a more subtle argument, which uses another result on unit equations, this time by Evertse, Schlickewei and Schmidt~\cite{ESS02}.  In addition, Evertse's theorem on $S$-unit equations does not apply to the rational functions $\Delta_n(z)$, and instead we work with related but slightly more complicated functions, see~\S3 for details\footnote{Note that there are some slight notational differences in ~\S3 between the present and published versions of this article.  These occur almost entirely with letters used for indices in various formul\ae.}.

\medskip
\subsection*{Context.}
Dynamical degrees play a key role in  algebraic, complex and arithmetic dynamics.
To any dominant rational map $f\colon X\dashrightarrow X$ of a projective variety
$X$ over $\bk$ is in fact associated a sequence $(\lambda_p(f))_{p=1}^{\dim X}$ of
dynamical degrees, each invariant under birational conjugation, see~\cite{DS05,Tru20,Dan20}; the dynamical degree above corresponds to $p=1$.  

Naturally defined in the context of algebraic dynamics, dynamical degrees were first introduced in \emph{complex} dynamics by Friedland~\cite{Fri91}, who showed that when $\bk=\bC$ and $f$ is a morphism, the topological entropy of $f$ is given by $\log\max_p\lambda_p(f)$; this generalized earlier work by Gromov, see~\cite{Gro03}, and was later extended (as an inequality) by Dinh and Sibony~\cite{DS05} to the case of dominant rational
maps.
Dynamical degrees are furthermore essential for defining and analyzing natural invariant currents and measures, see for example~\cite{RS97,Gue10,DS17} and the references therein.
Their importance from the point of view of complexity and integrability has also been exhibited in the physics literature by Bellon, Viallet and others, see~\eg~\cite{BV98,Via08}. 

In dimension two, the only relevant degrees are $\lambda_1 = \lambda$ and $\lambda_2$ (the `topological degree', equal to the number of preimages of a typical point if $\bk$ is algebraically closed of characteristic zero).  When $\bk=\bC$, their relationship determines which of two types of dynamical behavior (saddle or repelling) predominates (see~[DDG1-3] and~\cite{Gue05}).  The class of examples we consider here includes both types.  If, for instance, $\zeta = 1+2\imunit$, then we obtain a map $f$ of \emph{small topological degree} 
\begin{equation*}
\lambda_2(f) = \lambda_2(h) = |\zeta|^2 = 5 < \lambda_1(f) = 6.8575574092\dots
\end{equation*}
as computed numerically from Equation~($\star$).  Replacing $\zeta$ with $\zeta^2 = -3+4\imunit$, gives a map with \emph{large topological degree} $\lambda_2(f) = 25 > \lambda_1(f)=13.4496076817\dots$.

In \emph{arithmetic} dynamics, $\bk$ is a global field, and the (first) dynamical degree serves as an upper bound for the asymptotics of the growth of heights along orbits~\cite{Sil12,KS16,Mat20}; the question of when equality holds is part of the  \emph{Kawaguchi--Silverman conjecture}, which recently has attracted a lot of attention.

\medskip
\subsection*{Outlook}
As already mentioned, the set of all possible dynamical degrees is countable, and our Main Theorem shows that it contains transcendental numbers. It would obviously be interesting to say more about it. Note that the set of dynamical degrees of birational surface maps is much better understood, see \eg~\cite{BK06,McM07,Ueh16,BC16}.  It would be interesting to know (see \eg~\cite[page 1379]{Via08}) if a birational map $f\colon\bP^k\to\bP^k$ can have transcendental dynamical degree when $k\geq 3$.  
We intend to address this in a future article, though the number theoretic details seem more complicated. 
See~\cite{CaXi20,DF20} for some other results about degree growth of rational maps in higher dimensions.

\medskip
It would also be interesting to study the complex and arithmetic dynamics of the rational
map $f=f_\zeta$ considered here. For example, does $f$ admit a unique measure of maximal entropy, and is the topological entropy equal to $\log\lambda(f)$? The fact that $f$ is defined over $\bQ$ may be useful, see e.g.~\cite{JR18}, where it is shown that (complex) birational surface maps defined over $\bar\Q$ always admit a measure of maximal entropy $\log\lambda(f)$. 
On the arithmetic side, one may ask whether the Kawaguchi--Silverman conjecture holds: does every point with Zariski dense orbit have arithmetic degree equal to $\lambda(f)$? Note that what we call the Kawaguchi--Silverman conjecture is part~(d) of~\cite[Conjecture~6]{KS16}. Given our Main Theorem, the existence of a point as above would in fact contradict part~(b); see also~\cite{LS20}.

\medskip
\begin{ackn}
  We thank J. Blanc, S. Kawaguchi, C. T. McMullen, M. Satriano, U. Zannier and, especially, H. Krieger, for valuable comments.
  The last two authors thank J.-L. Lin and P. Reschke for their help during an
  earlier stage of this project, and J.~Lagarias and B.~Poonen for useful pointers
  regarding transcendence questions.
  Finally we thank the referee for a careful reading and many thoughtful suggestions.
  The first author was partially supported by NSERC grant RGPIN-2016-03632; the
  second author by NSF grant DMS-1954335;
  and the last author by NSF grants DMS-1600011 and DMS-1900025,
  and the United States---Israel Binational Science Foundation.
  The final form of the present collaboration originated at the Simons Symposium in Complex, Algebraic, and Arithmetic Dynamical Systems in May, 2019; we are very grateful to the Simons Foundation for its generous support.
\end{ackn}
%
%
%
%
\section{Dominant rational maps of the projective plane}
In this section we study dominant rational selfmaps of $\bP^2$ using the induced action on b-divisor classes.  The exposition largely follows~\cite{BFJ08, DL16}\footnote{Both of these articles were written for surfaces defined over $\C$, but the results we use from them work with proofs unchanged over any algebraically closed field.} but with particular attention paid to the structure of $\bP^2$ as a toric variety.  
We work over a field $\bk$ of characteristic different from two.  Since degrees of rational maps are invariant under ground field extension, we may and will assume that $\bk$ is algebraically closed. The assumption that $\charac\bk\ne2$ will be used in~\S\ref{S215}.

%
%
\subsection{Setup}.
Fix homogeneous coordinates $[x_0:x_1:x_2]$ on $\bP^2$ and use 
affine coordinates $(y_1,y_2)=(x_1/x_0,x_2/x_0)$
on the affine chart $\{x_0\ne0\}\simeq\bA^2$.
Recall that $\bP^2$ is a toric surface with torus $\mathbb{T}=\mathbb{G}^2_m=\{x_0x_1x_2\ne0\}$ and torus invariant prime divisors being the coordinate lines $\{x_j=0\}$, $j=0,1,2$. 
%
%
\subsection{Rational maps and their degrees}
A dominant rational selfmap of $\bP^2$ is given in homogeneous coordinates by
\begin{equation*}
  f\colon[x_0:x_1:x_2]\mapsto[f_0(x_0,x_1,x_2):f_1(x_0,x_1,x_2):f_2(x_0,x_1,x_2)],
\end{equation*}
where $f_0,f_1,f_2$ are homogeneous polynomials of the same degree $d\ge 1$, and with no factor in common.  The integer $\deg f:=d$ is called the \emph{degree} of $f$; see also Equation~\eqref{e301}.

The sequence $(\deg f^n)_{n\ge 1}$ is submultiplicative, \ie $\deg f^{m+n}\le\deg f^m\cdot\deg f^n$; hence the limit 
\begin{equation*}
\lambda(f)
=\lim_{n\to\infty}(\deg f^n)^{1/n}
=\inf_n(\deg f^n)^{1/n}\in[1,\infty)
\end{equation*} 
exists and
is equal to the dynamical degree of $f$ as defined in the introduction.
%
%
\subsection{Monomial maps}\label{S217}
Any $2\times 2$ matrix $\Linmap=(a_{ij})_{i,j}$ with integer coefficients and nonzero determinant defines a dominant rational self map $h_\Linmap\colon\bP^2\dashrightarrow\bP^2$, which in affine coordinates is given by 
$h_\Linmap\colon(y_1,y_2)\mapsto(y_1^{a_{11}}y_2^{a_{12}},y_1^{a_{21}}y_2^{a_{22}})$.
Such rational maps are called \emph{monomial maps}; they correspond to surjective endomorphisms of the algebraic group $\mathbb{T}$. 

Note that $h_{\Linmap_1\Linmap_2}=h_{\Linmap_1}\circ h_{\Linmap_2}$.  The degree of a monomial map is given by
\begin{equation*}
  \deg h_\Linmap
  =\max\{0,a_{11}+a_{12},a_{21}+a_{22}\}
  +\max\{0,-a_{11},-a_{12}\}
  +\max\{0,-a_{21},-a_{22}\},
\end{equation*}
see~\cite{HP07,BK08}.

\begin{figure}[ht]
  \centering
  \begin{tikzpicture}[scale=0.8]
    \draw (0,0) -- (3,0);
    \draw (0,0) -- (0,3);
    \draw (0,0) -- (0,-3);
    \draw (0,0) -- (-2.5,2.5);
    \draw (0,0) -- (-2.5,-2.5);
    \draw (2,1.5) node {$\re z+2\im z$};
    \draw (2,-1.5) node {$\re z-2\im z$};
    \draw (-1,2.1) node {$2\im z$};
    \draw (-1,-2.1) node {$-2\im z$};
    \draw (-2,0) node {$-2\re z$};
  \end{tikzpicture}
  \hspace*{1cm}
  \begin{tikzpicture}[scale=0.8]
    \draw (0,0) -- (3,0);
    \draw (0,0) -- (0,3);
    \draw (0,0) -- (0,-3);
    \draw (0,0) -- (-2.5,2.5);
    \draw (0,0) -- (-2.5,-2.5);
    \draw (2,1.5) node {$\gamma=1-2\imunit$};
    \draw (2,-1.5) node {$\gamma=1+2\imunit$};
    \draw (-1,2.1) node {$\gamma=-2\imunit$};
    \draw (-1,-2.1) node {$\gamma=2\imunit$};
    \draw (-2,0) node {$\gamma=-2$};
  \end{tikzpicture}
  \caption{The left picture shows the piecewise $\bR$-linear function $\psi\colon\C\to\R_{\geq 0}$ defined by Equation~\eqref{e208}. The right picture shows the element $\gamma\in\Gamma_0$ that realizes the maximum in the definition of $\psi$. Note that the angles that the rays make with the positive real axis are all integer multiples of $\pi/4$.}\label{F210}
\end{figure}
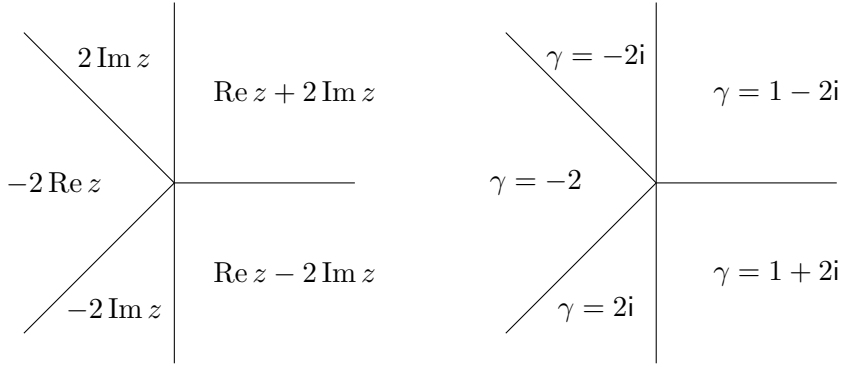

We can view $\Lambda$ as a linear selfmap of $\Z^2$ or $\R^2$.
Now identify $\R^2$ with $\C$ and assume that $\Lambda$ is given by multiplication with a Gaussian integer $\zeta\in\Z[\imunit]$, that is,
\begin{equation*}
  \Linmap=
  \Linmap_\zeta=
  \begin{pmatrix}
    \re\zeta & -\im\zeta\\
    \im\zeta & \re\zeta\\
  \end{pmatrix},
\end{equation*}
In this case, we write $h_\zeta=h_{\Linmap_\zeta}$. Note that $h_{\zeta_1}\circ h_{\zeta_2}=h_{\zeta_1\zeta_2}$. We have 
\begin{multline*}\label{e125}
  \deg h_\zeta
  =\max\{0,\re\zeta-\im\zeta,\re\zeta+\im\zeta\}\\
  +\max\{0,-\re\zeta,\im\zeta\}
  +\max\{0,-\re\zeta,-\Im\zeta\}.
\end{multline*}
which we can rewrite as $\deg h_\zeta=\psi(\zeta)$, where $\psi\colon\C\to\R_{\ge0}$ is a convex piecewise $\R$-linear function given by
\begin{equation}\label{e208}
  \psi(z):=\max_{\gamma\in\Gamma_0}\re(\gamma z),
\quad\text{where}\quad
  \Gamma_0:=\{-2,\pm2\imunit,1\pm2\imunit\};
\end{equation}
see Figure~\ref{F210}.

One checks that $\psi$ is comparable to the Euclidean norm on $\C$; specifically, $|z|\le\psi(z)\le \sqrt{5}|z|$.
Since $h_\zeta^n=h_{\zeta^n}$ for $n\ge1$, it follows that the
dynamical degree of $h_\zeta$ is
\begin{equation*}
  \lim_{n\to\infty}\psi(\zeta^n)^{1/n}=|\zeta|.
\end{equation*}

We will be interested in the case when $\zeta^n\not\in\R$ for all $n\ge1$. This is equivalent to $\zeta$ not being an integer multiple of $1$, $\imunit$ or $1\pm\imunit$, see \eg~\cite[Main Lemma]{Cal09}.  In this
case, there is, for every $n \geq 1$, a unique element $\gamma(n)\in \Gamma_0$ such that $\psi(\zeta^n) =
\re(\gamma(n)\zeta^n)$
%
%
\subsection{Blowups}
By a \emph{blowup} of $\bP^2$ we mean a birational morphism $\pi\colon X_\pi\to\bP^2$, where $X_\pi$ is a smooth projective surface.  Up to isomorphism, $\pi$ is then a finite composition of point blowups~\cite[Theorem~4.10]{Shafarevich}.  If $\pi$ and $\pi'$ are blowups of $\bP^2$, then $\mu:=\pi^{-1}\circ\pi'\colon X_{\pi'}\dashrightarrow X_\pi$ is a birational map; we say that $\pi'$ \emph{dominates} $\pi$, written $\pi'\ge\pi$, if $\mu$ is a morphism.  Any two blowups, can be dominated by a third, as follows by applying~\cite[Theorem~4.9]{Shafarevich} to the birational map $\mu$ above.  
It follows that the set $\mathfrak{Bl}$ of isomorphism classes of blowups is a directed set.
%
%
\subsection{Primes over $\bP^2$}  
We will say that prime divisors $E\subset X_\pi$ and $E'\subset X_{\pi'}$ in different blowups are \emph{equivalent} if there is a blowup $\pi'' = \pi\circ\mu = \pi'\circ\mu'$ dominating both $\pi$ and $\pi'$ and a prime divisor $E'' \subset X_{\pi''}$ such that  $E = \mu(E'')$ and $E' = \mu'(E'')$.  We let $\cP$ denote the set of all the resulting equivalence classes and call each  $E\in\cP$ a \emph{prime over $\bP^2$}.

We say that a blowup $\pi$ \emph{expresses} a prime $E\in\cP$, if $E$ is represented by a prime divisor in $X_\pi$ that we call then the \emph{center} of $E$ on $X_\pi$.  Slightly abusively, we use the same letter to denote the center, writing $E\subset X_\pi$.

If a blowup $\pi$ does not express a prime $E$, then we can choose a further blowup $\pi' = \pi\circ\mu>\pi$ such that $E$ is represented by a prime divisor on $X_{\pi'}$.  The image under $\mu$ of this prime divisor is a point in $X_\pi$ which does not depend on the choice of $\pi'$ and which we call the center of $E$ on $X_\pi$.
%
%
\subsection{b-divisor classes}
For any blowup $\pi$ of $\bP^2$, denote by $\Pic(X_\pi)$ the Picard group on $X_\pi$, \ie the set of linear equivalence classes of (Cartier) divisors on $X_\pi$.  When $\pi'\ge\pi$, the birational morphism $\mu\colon X_{\pi'}\to X_\pi$ induces an injective homomorphism $\mu^*\colon\Pic(X_{\pi})\to\Pic(X_{\pi'})$.  The group of \emph{b-divisor classes} on $\bP^2$ is defined as the direct limit
\begin{equation*}
  \cC:=\varinjlim_{\pi\in\mathfrak{Bl}}\Pic(X_\pi).
\end{equation*}
Concretely, an element of $\cC$ is an element of $\Pic(X_\pi)$ for some blowup $\pi$,
where two elements $A\in\Pic(X_\pi)$, $A'\in\Pic(X_{\pi'})$ are identified iff they pull back to the same class on some blowup dominating both $\pi$ and $\pi'$.
A class in the image of $\Pic(X_{\pi})\hookrightarrow\cC$ is said to be \emph{determined} on $X_\pi$.  We let
\begin{equation*}
  H=\cO_{\P^2}(1)\in\cC
\end{equation*}
denote the class determined by a line in $\bP^2$.

\begin{remark}
The $b$ in `$b$-divisor' stands for birational, following Shokurov. In~\cite{BFJ08}, the elements of $\cC$ were referred to as Cartier classes on the Riemann--Zariski space of $\bP^2$. The space $\cC$ appears earlier in~\cite{Manin}, where it is denoted $Z^\bullet(\P^2)$. Note that since each surface $X_\pi$ is rational, the Picard group $\Pic(X_\pi)$ coincides with the N\'eron-Severi group $\NS(X_\pi)$.
\end{remark}

There is a natural intersection pairing $\cC\times\cC\to\bZ$, denoted $(A\cdot B)$ for $A,B\in\cC$. This is defined as the intersection number on any blowup $X_\pi$ where $A$ and $B$ are both determined
(see~\cite[\S1.4]{BFJ08} or~\cite[\S34.7]{Manin}).
%
%
\subsection{Toric blowups}
We call a blowup $\pi\colon X_\pi\to\bP^2$ \emph{toric} if $X_\pi$ is also a toric surface and $\pi$ is equivariant with respect to the torus action.  Concretely, $\pi = \pi_1\circ\dots\circ\pi_n$, where each factor $\pi_j\colon X_j\to X_{j-1}$ is a point blowup centered at the intersection of two different torus invariant prime divisors in $X_{j-1}$.

If $\pi$ is a toric blowup of $\P^2$ and $E\subset X_\pi$ a torus invariant prime divisor, then a point $p\in E$ is called \emph{free} if it does not belong to any other torus invariant prime divisor on $X_\pi$, \ie its orbit under the torus action is 1-dimensional.

We will call $E\in\cP$ a \emph{toric prime} if there is a toric blowup $\pi$ that expresses $E$ as a torus invariant prime divisor.  Let $\cPtor$ denote the set of all toric primes.  

\begin{proposition}
\label{prop:factor}
Any blowup $\pi$ of $\bP^2$ factors uniquely as $\pi = \pitor\circ\mu$ into a toric blowup $\pitor$ that expresses the same set of toric primes as $\pi$ and a birational morphism $\mu\colon X_\pi \to X_{\pitor}$ that contracts only non-toric primes.  
\end{proposition}

\begin{proof}
This follows from~\cite[Corollary 5.5]{DL16} and the fact that the toric primes in $X_\pi$ are precisely the (simple) poles of the rational $2$-form $\pi^*\frac{dy_1\wedge dy_2}{y_1y_2}$.
\end{proof}

To each prime $E\in\cP$ we associate an order of vanishing valuation $\ord_E\colon\bk(\bP^2)^\times\to \Z$ by choosing a blowup $\pi$ such that $E\subset X_\pi$ and setting $\ord_E(\varphi)$ equal to the coefficient of $E$ in the divisor of the rational function $\varphi\circ\pi$ on $X_\pi$.  We define a `tropicalization' map $\trop\colon\cP\to \Z^2$ by
\begin{equation*}
\trop(E) = (\ord_E(y_1),\ord_E(y_2)),
\end{equation*}
where $(y_1,y_2)$ are the affine coordinates fixed above.  Note that $\trop(E)=(0,0)$ for all non-toric primes $E\subset\bP^2$, whereas if we write $H_j:=\{x_j=0\}$, $j=0,1,2$, then 
\begin{equation}\label{e302}
  \trop(H_0) = (-1,-1),\quad \trop(H_1) = (1,0), \quad \trop(H_2)=(0,1).
\end{equation} 
In any blowup $\pi$ of $\bP^2$, the divisors of the rational functions $y_j\circ\pi$, $j=1,2$, have simple normal crossings support.  Hence, if $\mu\colon X_{\pi'}\to X_{\pi}$ is the point blowup at $p\in X_{\pi}$, then the prime $E'$ contracted by $\mu$ satisfies
\begin{equation}
\label{eqn:blowuptrop}
\trop(E') = \sum_{E\subset X_\pi\, \text{s.t.\ $p\in E$}} \trop(E).
\end{equation}
We call a non-zero element $t\in\Z^2$ \emph{primitive} if $t\not\in m\Z^2$ for any integer $m\ge2$. The next result follows easily by induction from the discussion above and is related to the fact that $SL_2(\bZ)$ acts transitively on primitive elements of $\Z^2$.
\begin{proposition}
\label{prop:primitive}
The map $\trop$ restricts to a bijection from $\cPtor$ onto the set of primitive elements $t\in\Z^2$.  
\end{proposition}
We will say that elements $s,t\in\Z^2$ are \emph{commensurate} if $s = rt$ for some positive $r\in\bQ$.  For each non-zero (but not necessarily primitive) element $t\in\Z^2$, we let $E_t\in\cPtor$ be the unique toric prime such that $\trop(E_t)$ is commensurate with $t$.
\begin{proposition}
\label{prop:trope}
Let $\pi = \pitor\circ\mu$ be a blowup of $\P^2$, factored as in Proposition~\ref{prop:factor}, and $E\subset X_\pi$ a non-toric prime divisor with $t:=\trop(E)\neq (0,0)$.  Then
\begin{enumerate}
\item $\pitor$ expresses the toric prime $E_t$, and $\mu(E)$ is a free point on $E_t$;
\item if $E'\in\cP$ is not expressed in $X_\pi$, and its center on $X_\pi$ is a point $p'\in E$, then $E'$ is also non-toric, and $\trop(E')$ is commensurate with $t$.
\end{enumerate}
\end{proposition}
\begin{proof}
The support of the divisor of the rational function $y_j\circ\pitor$ on $X_{\pitor}$ does not meet the torus $\mathbb{T}$ for $j=1,2$, so since $E\notin\cPtor$ and $\trop(E)\ne(0,0)$, we have that $\mu(E)$ is a point in $X_{\pitor}\setminus\mathbb{T}$. If $p$ is the intersection of two distinct toric primes expressed by $\pitor$, then $\pi$ dominates $\pitor\circ\mu'$ where $\mu'\colon X_{\pi'}\to X_{\pitor}$ is the point blowup at $p$.  This means, however, that ${\mu'}^{-1}(p)$ is a toric prime expressed by $\pi$ but not $\pitor$, which contradicts the choice of $\pitor$.

Thus $p$ is a free point on a prime $E_s$ expressed by $\pitor$, with $s=\trop(E_s)$.  The map $\trop$ therefore vanishes along all other primes expressed by $\pitor$ that contain $p$.  Hence by factoring $\mu$ into point blowups and repeatedly applying Equation~\eqref{eqn:blowuptrop}, we see that $t=\trop(E)$ is commensurate with $s$. So~(i) holds, and we turn to (ii).

By the previous step, any prime that is expressed by $\pi$ and contains $p'$ has tropicalization equal to a multiple (possibly $0$) of $t$.  Hence we can choose a blowup $\pi' = \pi\circ\mu'>\pi$ that expresses $E'$, factor $\mu'$ into point blowups and repeatedly apply Equation~\eqref{eqn:blowuptrop} to obtain that $\trop(E')$ is commensurate with $t$.  Since $E'\neq E_t$, Proposition~\ref{prop:primitive} tells us that $E'$ is not toric.
\end{proof}

The set of \emph{toric b-divisor classes} $\cCtor\subset\cC$ is the direct limit $\varinjlim_\pi\Pic(X_\pi)$, where $\pi$ runs over all toric blowups of $\bP^2$. Each class in $\Pic(X_\pi)$ is represented by a toric divisor, \ie a divisor with support equal to a collection of toric primes expressed by $\pi$.  In particular $H\in\cCtor$, and a class in $\cCtor$ is orthogonal to $H$ iff it is represented by a $\pi$-exceptional toric divisor on some toric blowup $\pi$ of $X$.  We will use this fact below in proving Lemma~\ref{L104}.
%
%
\subsection{Action by rational maps on primes and on b-divisor classes}
Consider a dominant rational map $f\colon\bP^2\dashrightarrow\bP^2$. For any blowups $\pi,\pi'$ of $\bP^2$ we have an induced rational map
$f_{\pi\pi'}:=\pi^{-1}\circ f\circ\pi'\colon X_{\pi'}\dashrightarrow X_\pi$.
Given $\pi$, we can choose $\pi'$ such that $f_{\pi\pi'}$ is a morphism, as follows from~\cite[Theorem~4.8]{Shafarevich}.
We now define a group homomorphism
\begin{equation*}
  f^*\colon\cC\to\cC
\end{equation*}  
as follows: if $A\in\cC$ is determined on $X_\pi$, pick a blowup $\pi'$ such that $f_{\pi\pi'}\colon X_{\pi'}\to X_\pi$ is a morphism, and declare $f^*A\in\cC$ to be the class determined on $X_{\pi'}$ by $f_{\pi\pi'}^*A$.
This action is functorial: if $f$ and $g$ are dominant rational 
maps of $\bP^2$, then $(f\circ g)^*=g^*f^*$ on $\cC$.
When $h\colon\bP^2\dashrightarrow\bP^2$ is monomial, we have  $h^*\cCtor\subset\cCtor$.
The degree of a rational map can be computed as follows:
\begin{equation}\label{e301}
  \deg f = (f^*H\cdot H).
\end{equation}

The rational map $f$ also induces an action $f\colon\cP\to \cP$ on the set of all primes over $\bP^2$.  If $\pi'$ is a blowup expressing $E\in \cP$, then as in~\cite{BFJ08} (see just before Lemma 2.4) there exists  another blowup $\pi$ such that the lift $f_{\pi\pi'}\colon X_{\pi'}\tto X_{\pi}$ does not contract any curves.  We set $f(E) := f_{\pi\pi'}(E)$. 

\begin{proposition}
\label{prop:monomialaction}
For any monomial map $h\colon\bP^2\tto \bP^2$, with associated matrix $\Lambda$, and any prime $E\in \cP$, 
\begin{enumerate}
 \item $h(E)$ is toric if and only if $E$ is; and
 \item $\trop(h(E))$ is commensurate with $\Lambda(\trop(E))$.
\end{enumerate}
\end{proposition}

\begin{proof}
  The first conclusion follows from the first conclusion of~\cite[Corollary 6.3]{DL16} and the fact that $h^* \frac{dy_1\wedge dy_2}{y_1y_2} = (\det\Lambda)\,\frac{dy_1\wedge dy_2}{y_1y_2}$.  The second conclusion is a (by now) standard computation.
\end{proof}
%
%
%
%
\section{The degree sequence of certain rational maps}\label{S215}
We now specialize the considerations above to a particular class of maps that will later be shown to have transcendental dynamical degrees.
%
%
\subsection{A volume preserving involution}
As in~\cite{DL16} we consider the involution\footnote{Here we use that the ground field has characteristic different from two. Indeed, $g$ is the identity in characteristic two.} $g\colon\bP^2\dashrightarrow\bP^2$
defined in homogeneous coordinates by
\begin{equation*}
g\colon[x_0:x_1:x_2]\mapsto[x_0(x_1+x_2-x_0):x_1(x_2+x_0-x_1):x_2(x_0+x_1-x_2)].
\end{equation*}
In affine coordinates $(y_1,y_2)=(x_1/x_0,x_2/x_0)$, this becomes
\begin{equation}\label{e220}
  g\colon
  (y_1,y_2)\mapsto \left(-y_1\frac{1-y_1+y_2}{1-y_1-y_2},-y_2\frac{1+y_1-y_2}{1-y_1-y_2}\right).
\end{equation}

The projective linear automorphism
\begin{equation*}
  A\colon[x_0:x_1:x_2]\mapsto[x_1+x_2-x_0:x_2+x_0-x_1:x_0+x_1-x_2]
\end{equation*}
conjugates $g$ to the Cremona involution $AgA^{-1}\colon[x_0:x_1:x_2]\mapsto[x_1x_2:x_2x_0:x_0x_1]$. As a consequence, we have the following geometric description. Consider the three points $p_0=[0:1:1]$, $p_1=[1:0:1]$, $p_2=[1:1:0]$ and the three lines $L_0=\{x_0=x_1+x_2\}$, $L_1=\{x_1=x_2+x_0\}$,  $L_2=\{x_2=x_0+x_1\}$ on $\bP^2$. Let $X^0$ be the blowup of 
$\bP^2$ at $\{p_0,p_1,p_2\}$, with exceptional divisors $F_0,F_1,F_2$. Then $g$ induces an automorphism of $X^0$ of order two that sends $F_j$ to the strict transform of $L_j$ for $j=0,1,2$.
\begin{figure}[ht]
  \centering
  \begin{tikzpicture}[scale=1.0]
    \tikzmath{\dotsize=1.0; }

    \draw (-3,2) -- (3,2);
    \draw (-3,3) -- (1,-1);
    \draw (3,3) -- (-1,-1);

    \draw[blue] (1,5) -- (-5,-1);
    \draw[blue] (-1,5) -- (5,-1);
    \draw[blue] (-5.5,0) -- (5.5,0);

    \node at (-2,2) [circle,fill,inner sep=\dotsize]{};
    \node at (2,2) [circle,fill,inner sep=\dotsize]{};
    \node at (0,0) [circle,fill,inner sep=\dotsize]{};

    \node at (0,4) [circle,fill,inner sep=\dotsize]{};
    \node at (-4,0) [circle,fill,inner sep=\dotsize]{};
    \node at (4,0) [circle,fill,inner sep=\dotsize]{};

    \draw (3.5,2) node {$L_0$};
    \draw (-3.5,3) node {$L_2$};
    \draw (-1.5,-1) node {$L_1$};

    \draw (-1.2,4.5) node {$H_2$};
    \draw (-4.5,-1) node {$H_1$};
    \draw (5.5,0.2) node {$H_0$};

    \draw (0,-0.5) node {$p_0$};
    \draw (-2,1.5) node {$p_1$};
    \draw (2,2.5) node {$p_2$};

  \end{tikzpicture}
  \caption{The birational involution $g$ contracts the line $L_j$ to the point $p_j$, $j=0,1,2$. It leaves the coordinate lines $H_j=\{x_j=0\}$ invariant. The restriction $g|_{H_j}$ fixes the two points $H_j\cap H_l$, $l\ne j$, and sends $p_j$ to the point $H_j\cap L_j$ (which is not shown).}\label{F204}
\end{figure}
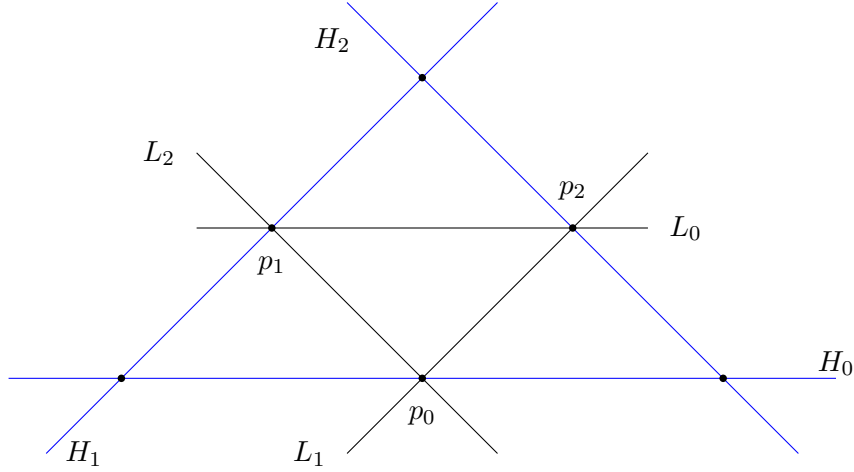 

Let $\pi:X_\pi\to\cP^2$ be a toric blowup.  For $i=1,2,3$, the point $p_i\in \bP^2$ is a free point on the toric prime $H_i$.  Hence its preimage by $\pi$ remains a free point on $H_i$, and we continue to denote it by $p_i$.  We let $\pi^0\colon X^0_\pi\to\bP^2$ be the blowup of $X_\pi$ along the set $\{p_0,p_1,p_2\}\subset X_\pi$.

\begin{lemma}\label{L109}
For any toric blowup $\pi\colon X_\pi\to\bP^2$, the induced birational map $g^0_{\pi\pi}\colon X^0_\pi\dashrightarrow X^0_\pi$ is a morphism that fixes each toric prime $E\subset X^0$.
\end{lemma}

\begin{proof}
We have already explained that this is true when $X_\pi = \bP^2$.  Hence it suffices by induction to show that if the lemma holds for some toric blowup $\pi$, then it also holds for the toric blowup $\pi' = \pi\circ\mu$, where $\mu$ is the point blowup of the intersection of two toric primes $E_s,E_t\subset X_\pi$.  
But the facts that $g$ fixes both $E_s$ and $E_t$ and that $E_s\cap E_t$ is distinct from $p_1,p_2,p_3$ imply that $E_s\cap E_t$ remains a point in $X_\pi^0$ and that the automorphism $g^0_{\pi\pi}$ fixes it.  Hence $g^0_{{\pi'}{\pi'}}$ is an automorphism fixing the exceptional prime $E = \mu^{-1}(E_s\cap E_t)$.  
\end{proof}

\begin{lemma}\label{L102} The induced map
$g\colon\cP\to \cP$ is a bijection that fixes the subset $\cPtor$ pointwise.  If $E\in\cP$ is a prime such that $\trop(E)$ is incommensurate with $(0,0)$, $(-1,-1)$, $(1,0)$, and $(0,1)$, then $\trop(g(E))$ is commensurate with $\trop(E)$.
\end{lemma}
Though it is not strictly necessary for the proof we note the related fact that $Dg(q_{ij}) = -I$ at each fixed point $q_{ij} = H_i\cap H_j$. 

\begin{proof}
Set $t=\trop(E)$.  The first assertion follows from Lemma~\ref{L109}.  For the second assertion, we may assume that $E$ is not toric.  If $\pi = \pitor\circ\mu$ is a blowup that expresses $E$, decomposed as in Proposition~\ref{prop:factor}, then Proposition~\ref{prop:trope} tells us that $\mu(E)$ is a free point on the toric prime $E_t \subset X_{\pitor}$.  Since $t$ is incommensurate with $(-1,-1)$, $(1,0)$ and $(0,1)$, we have $\mu(E)\notin \{p_1,p_2,p_3\}$; see Equation~\eqref{e302}. So by Lemma~\ref{L109}, the map $g_{\pitor\pitor}$ is a local isomorphism about $\mu(E)$ and the image $g_{\pitor\pitor}(\mu(E))$ is also a free point in $E_t$.  Thus $g(E)$ is a non-toric prime over a free point in $E_t$, and $\trop(g(E))$ is commensurate with $t$.
\end{proof}

Now consider a monomial map $h=h_\zeta$ associated (as in~\S\ref{S217}) to a Gaussian integer $\zeta$  for which $\zeta^n\not\in\R$ for all $n\ge1$. We will construct a set $\cP'\subset \cP$ of primes over $\bP^2$ that is backward invariant under both $g$ and $h$.  As before, we identify $\R^2$ with $\C$.
Define
\begin{equation*}
  N'
  :=\bigcup_{j\ge1}\zeta^{-j}\R_{>0} \{1,\imunit,-1-\imunit\}
  \subset\C.
\end{equation*}
Our assumption on $\zeta$ implies that $N'$ is an infinite set of rays in $\C$, none containing $0$, $1$, $\imunit$, or $-1-\imunit$.  Let
\begin{equation*}
  \cP':= \{E\in \cP\setminus \cPtor\mid \trop(E) \in N'\}.
\end{equation*}

\begin{corollary}\label{C105}
We have $g^{-1}(\cP')\subset \cP'$ and $h^{-1}(\cP')\subset \cP'$.
\end{corollary}

\begin{proof}
The first inclusion follows from Lemma~\ref{L102} and the fact that $1,\imunit,-1-\imunit\notin N'$.  The second inclusion follows from Proposition~\ref{prop:monomialaction} and the fact that the matrix $\Lambda$ associated to the monomial map $h$ acts on $\R^2\simeq\C$ by multiplication with $\zeta$.
\end{proof}

Next we study the action of $g$ and $h$ on the group $\cC$ of b-divisor classes.  Define $\cC'\subset\cC$ to be the subgroup of classes that can be represented by a divisor $D$ on some blowup of $\bP^2$, such that all irreducible components of $D$ lie in $\cP'$. Proposition~\ref{prop:trope} implies that $\cC'$ is orthogonal to $\cCtor$ and in particular to $H=\cO_{\P^2}(1)$.
\begin{corollary}\label{C108}
  We have $g^*\cC'\subset\cC'$ and $h^*\cC'\subset\cC'$.
\end{corollary}

\begin{proof}
By linearity it suffices to consider the pullback of a prime divisor $E\subset X_\pi$ with $\trop(E)\in N'$.  If $\pi'$ is a blowup of $\P^2$ such that $g_{\pi\pi'}\colon X_{\pi'}\tto X_\pi$ is a morphism, then $g^*E$ is determined in $X_\pi'$ by $g_{\pi\pi'}^*E$.  Further, every irreducible component $E'$ of $g_{\pi\pi'}^*E$ satisfies $g_{\pi\pi'}(E') \subset E$.  Thus, as elements of $\cP$, either $g(E') = E$ or the center of $g(E')$ on $X_\pi$ is a point in $E$.  In the second case, Proposition~\ref{prop:trope} implies that $g(E')$ is non-toric with $\trop(g(E'))\in N'$.  Hence, in either case, Corollary~\ref{C105} yields $E'\in\cP'$ and therefore $g_{\pi\pi'}^*E\in \cC'$.  The proof that $h^*\cC'\subset\cC'$ is identical.
\end{proof}

Next we study the action of $g$ on toric b-divisor classes.  

\begin{lemma}\label{L103}
  We have $g^*H=2H+R$, where $R\in\cC$ satisfies $h^*R\in\cC'$.
\end{lemma}

\begin{proof}
We use the notation introduced earlier in the subsection. On $X^0$, $H$ is represented by the divisor $\frac13\sum_{i=0}^2(L_i+2F_i)$, so $g^*H$ is represented by $\tfrac13\sum_{i=0}^2(2L_i+F_i)=2H+R$, where $R:=-\sum_{i=0}^2F_i$. It only remains to see that $h^*R\in\cC'$.  Pick a blowup $\pi\colon X_\pi\to\bP^2$ such that $h$ induces a morphism $h_\pi\colon X_\pi\to X^0$. Then $h^*R$ is represented by the divisor $\sum_{i=0}^2h_\pi^*F_i$ on $X_\pi$.  Every irreducible component $F\subset X_\pi$ of $h_\pi^* F_i$ satisfies  $h_\pi(F) \subset F_i$.  Applying Proposition~\ref{prop:trope} if $h_\pi(F)$ is a point, we find that the prime $h(F)\in\cP$ is non-toric, and $\trop(h(F))$ is commensurate with $-1-\imunit$, $1$, or $\imunit$.  Proposition~\ref{prop:monomialaction} implies then that $F$ is a non-toric prime with $\zeta \trop(h(F))$ commensurate with $-1-i$, $i$, or $1$.  We conclude that $F\in \cP'$, and $h^* F_i\in \cC'$.
\end{proof}

\begin{lemma}\label{L104}
  If $A\in\cCtor$ and $(A\cdot H)=0$, then $g^*A=A$.
\end{lemma}

\begin{proof}
  There exists a toric blowup $\pi\colon X_\pi\to\bP^2$ such that $A$ is represented by a torus invariant $\pi$-exceptional divisor on $X_\pi$.  Let $\mu\colon X^0_\pi\to X_\pi$ be the blowup of $X_\pi$ at $p_0,p_1,p_2$.  Since no irreducible component of $A$ in $X_\pi$ is the proper transform of one of the coordinate lines $H_i$, it follows that $\mu^*A$ is still supported on toric primes in $X^0_\pi$. By Lemma~\ref{L109}, the birational map $g^0_{\pi\pi}\colon X^0_\pi\dashrightarrow X^0_\pi$ is a morphism, and $(g^0_{\pi\pi})^*\mu^*A=\mu^*A$ in $\Pic(X^0_\pi)$. This implies $g^*A=A$ in $\cC$.
\end{proof}

%
%
\subsection{Degree sequence}
Let $g$ be the involution above, $h=h_\zeta$ the
monomial map associated to a Gaussian integer $\zeta$ such that
$\zeta^n\not\in\bR$ for all $n\ge1$, and set
\begin{equation*}
  f:=g\circ h.
\end{equation*}
Write
\begin{equation*}
  d_n=\deg(h^n)=(h^{n*}H\cdot H)
  \quad\text{and}\quad
  e_n=\deg(f^n)=(f^{n*}H\cdot H)
\end{equation*}
for $n\ge 0$. In particular, $d_0=e_0=1$.
Our aim is to prove the following recursion formula.
\begin{proposition}\label{P106}
  We have $e_n=d_n+\sum_{j=0}^{n-1}e_jd_{n-j}$ for $n\ge 0$.
\end{proposition}
\begin{proof}
  We will prove the following more precise result by induction on $n$.
  \begin{align}
    f^{n*}H
    &=h^{n*}H+\sum_{j=0}^{n-1}e_jh^{(n-j)*}H\
      \bmod \cC'\tag{$A_n$}\\
    g^*f^{n*}H
    &=h^{n*}H+\sum_{j=0}^ne_jh^{(n-j)*}H+e_nR\ 
      \bmod \cC'\tag{$B_n$}
  \end{align}
  Pairing~($A_n$) with $H$ implies the desired result since $\cC'$ is orthogonal to $H$.

  Now $(A_0)$ is trivial, and $(B_n)\Rightarrow(A_{n+1})$ for $n\ge 0$, as is seen by applying $h^*$ and using that $h^*R\in\cC'$. It therefore suffices to prove that   $(A_n)\Rightarrow(B_n)$ for $n\ge0$.

  To this end, we rewrite $(A_n)$ as
  \begin{equation*}
    f^{n*}H
    =e_nH+\left(h^{n*}H+\sum_{j=0}^{n-1}e_jh^{(n-j)*}H-e_nH\right)
    \bmod\cC'.
  \end{equation*}
  The expression in parentheses lies in $\cCtor$ and is orthogonal
  to $H$. Lemmas~\ref{L103} and~\ref{L104} therefore give
  \begin{align*}
    g^*f^{n*}H
    &=2e_nH+e_nR
    +(h^{n*}H+\sum_{j=0}^{n-1}e_jh^{(n-j)*}H-e_nH)\
    \bmod\cC'\\
    &=h^{n*}H+\sum_{j=0}^ne_jh^{(n-j)*}H+e_nR\
      \bmod\cC',
  \end{align*}
  which completes the proof.  
\end{proof}
%
%
\subsection{Dynamical degree}
Set $\Delta_h(z):=\sum_{j=1}^\infty d_jz^j$ and $\Delta_f(z):=\sum_{j=1}^\infty e_jz^j$.
These are power series with radii of convergence equal to $|\zeta|^{-1}$ and $\lambda^{-1}$, respectively, where $\lambda$ is the dynamical degree of $f$. 
Proposition~\ref{P106} shows that
\begin{equation}\label{e214}
  (2+\Delta_f(z))(1-\Delta_h(z))=2
\end{equation}
for $|z|<\min\{\lambda^{-1},|\zeta|^{-1}\}$.
\begin{proposition}\label{P107}
  The dynamical degree $\lambda=\lambda(f)$ satisfies $\lambda>|\zeta|$,
  and $\lambda$ is the unique positive solution to
  the equation $\sum_{j=1}^\infty d_j\lambda^{-j}=1$, where $d_j=\deg h^j$.
\end{proposition}
\begin{proof}
By submultiplicativity, $|\zeta| = \lim d_j^{1/j} = \inf d_j^{1/j}$; hence $d_j \geq |\zeta|^j$ for all $j$. Thus $\Delta_h(t)$ is positive and strictly increases from $0$ to $\infty$ on the interval $(0,|\zeta|^{-1})$.
Similarly, $\Delta_f(t)$ increases from $0$ to $\infty$ on $(0,|\lambda|^{-1})$.  The equation $2+\Delta_f(t)=\frac{2}{1-\Delta_h(t)}$ therefore implies that $t=|\lambda|^{-1}$ is the unique element of $(0,|\zeta|^{-1})$ for which $\Delta_h(t) = 1$.
\end{proof}

Now recall from~\S\ref{S217} that $d_j=\psi(\zeta^j)$, where $\psi$ is a convex, nonnegative  and piecewise $\bR$-linear function on $\C$ given by Equation~\eqref{e208} and illustrated in Figure~\ref{F210}.
Set
\begin{equation}\label{e212}
  \alpha=\lambda^{-1}\zeta.
\end{equation}
Then $|\alpha|<1$, and $\alpha$ is a solution to the equation
\begin{equation}\label{e213}
  1 = \re \Phi(\alpha),
\end{equation}
where $\Phi$ is a complex analytic function on the unit disk given by
\begin{equation}\label{e206}
  \Phi(z) := \sum_{j=1}^\infty \gamma(j) z^j,
\end{equation}
and where the coefficient $\gamma(j)$ is the element $\gamma\in\Gamma_0$ for which $\re(\gamma\alpha^j)$,
or equivalently $\re(\gamma\zeta^j)$, is maximized, see Figure~\ref{F210}.
If we write
\begin{equation*}
  \theta
  =\tfrac1{2\pi}\arg(\alpha)
  =\tfrac1{2\pi}\arg(\zeta)
  \in(0,1)\setminus\Q,
\end{equation*}
it follows that $\gamma(j)$ only depends on the image of $j\theta$ in $\bR\bmod\bZ$,
and more specifically which interval $(k/8,(k+1)/8)$ contains $j\theta\bmod 1$.
%
%
%
%
\section{Proof of transcendence}
We will spend the remainder of this article proving by contradiction that the number $\alpha$ in Equation~\eqref{e212}, and therefore the dynamical degree $\lambda(f)$, is transcendental. All that really matters going forward is that $|\alpha|<1$, that $\theta = \frac{\arg \alpha}{2\pi}$ is irrational, and that $\re \Phi(\alpha) = 1$, where $\Phi(z)$ is given by Equation~\eqref{e206}. Our arguments will be purely number theoretic, making no further use of algebraic geometry or dynamics. 
%
%
\subsection{Setup}\label{sec:setup} 
Since $\theta\notin\Q$, the sequence $(\gamma(j))_{j\ge 1}$ is aperiodic.  Nevertheless, as we will make precise below, it comes close to being $n$-periodic when $n$ is chosen to be the denominator in some continued fraction
approximant $m/n$ of $\theta$.  For such $n$, it will be illuminating to compare the analytic function $\Phi(z)$ with approximations by rational functions of the form
\begin{equation*}
  \Phi_n(z)
  := (1 - z^n)^{-1} \sum_{1\leq j \leq n} \gamma(j) z^j
  = \sum_{j\geq 1} \gamma_n(j) z^j,
\end{equation*}
where $\gamma_n(j)$ denotes the $n$-periodic extension of the initial sequence $\gamma(1),\dots,\gamma(n)$.
\begin{lemma}\label{L118}
For any sufficiently large $n\in\N$, we have $0<\re \Phi_n(\alpha)<1$.
\end{lemma}
\begin{proof}
By definition, we have 
\begin{equation*}
    1 - \re \Phi_n(\alpha)
    = \re(\Phi(\alpha) - \Phi_n(\alpha))
    = \sum_{j>n} \re((\gamma(j) - \gamma_n(j))\alpha^j).
\end{equation*}
Since $|\alpha|<1$ and $\Gamma_0=\{-2,\pm 2\imunit,1\pm 2\imunit\}$ is finite, the right side tends to zero as $n\to\infty$; in particular, $\re \Phi_n(\alpha)>0$ for large $n$.  Now, for each $j$, $\gamma(j)$ maximizes $\re(\gamma\alpha^j)$ over $\gamma\in\Gamma_0$, so $\re((\gamma(j) - \gamma_n(j))\alpha^j) \geq 0$. Thus $\re \Phi_n(\alpha)\le 1$, and to see that the inequality is strict, it suffices to find a single $j$ such that $\re((\gamma(j) - \gamma_n(j))\alpha^j) \neq 0$.  Since $\theta\notin\Q$, we can find $p\ge 1$ such that $p\theta\in(7/8,1)\bmod 1$. Assume $n>p$, and pick $m\ge 1$ such that if $j=mn+p$, then $j\theta\in(0,1/8)\bmod 1$. Then from Figure~\ref{F210}, we see that 
\begin{equation*}
\gamma(j) - \gamma_n(j) = \gamma(j)-\gamma(p) = (1-2\imunit)-(1+2\imunit)  = -4\imunit.
\end{equation*}
Since $\arg(\alpha^j)\in(0,\pi/4)$, it follows that $\re(-4\imunit\alpha^j)>0$.
\end{proof}

Lemma~\ref{L118} tells us $0< \re(\Phi(\alpha) - \Phi_n(\alpha)) < 1$ for large $n$.  To obtain better bounds, we clear the denominator in the definition of $\Phi_n(z)$, setting
\begin{align}
\label{eqn:error}
  \Psi_n(z)
  &:= 2|1-z^n|^2 \re (\Phi(z) - \Phi_n(z))\notag\\
  &= 2\Re\left((1-\bar z^n)\left(
   \sum_{j=1}^\infty (1-z^n)\gamma(j) z^j - \sum_{j=1}^n \gamma(j) z^j 
   \right)\right)\notag\\
  &= 2\Re\left((1-\bar z^n)\left(
   \sum_{j=n+1}^\infty \gamma(j) z^j - z^n\sum_{j=1}^\infty \gamma(j) z^j 
   \right)\right)\notag\\
  &= 2\Re\left((1-\bar z^n)\sum_{j = n+1}^\infty (\gamma(j) - \gamma(j-n)) z^j\right).
\end{align}
Since $|\alpha|<1$, we have that $0<\Psi_n(\alpha) <1$ for large $n$.  The final expression for $\Psi_n$ makes the following terminology convenient.
\begin{definition}
  We say that an index $j > n$ is \emph{$n$-regular} if $\gamma(j) = \gamma(j-n)$, and \emph{$n$-irregular} otherwise.  
\end{definition}

Since $\theta$ is irrational, there are infinitely many $n$-irregular indices, but they 
nevertheless form a rather sparse subset of $\bN$, as will be explored below.  Our arguments will depend on how well $\theta$ can be approximated by rational numbers.  Recall (from e.g. Chapters X-XI of~\cite{HW}) that any irrational number $t\in\R$ admits an infinite sequence of continued fraction approximants $m_i/n_i$, with $n_i$ strictly increasing, $m_i$ coprime to $n_i$, and $|n_i t - m_i|< \frac 1{n_i}$ for all $i\in\N$. 
\begin{proposition}
\label{prop:bad}
Let $t\in\R$ be irrational with continued fraction approximants $m_i/n_i$, $i\in\N$.  Then the following are equivalent.
\begin{enumerate}
 \item There exists $\kappa>0$ such that $|n_i t - m_i| \geq \frac{\kappa}{n_i}$ for all $i\in\N$.
 \item There exists $\kappa>0$ such that $|n t - m| \geq \frac{\kappa}{n}$ for all $m,n\in\Z$ with $n>0$.
 \item There exists $A$ such that $n_{i+1} \leq An_i$ for all $i\in\N$.
 \item The coefficients in the continued fraction expansion of $t$ are uniformly bounded.
\end{enumerate}
\end{proposition}
\begin{proof}
Suppose first that (i) holds.  For each $i$ we have
$|n_i t-m_i| < 1/n_{i+1}$~\cite[Corollary 1.4]{Bugeaud}, and hence $\kappa/n_i < 1/n_{i+1}$, which gives that (iii) holds with $A=\kappa^{-1}$.  Next suppose that (iii) holds. If $a_j$ is the $j$-th coefficient in the continued fraction of $t$ then $n_{i+1} = a_{i+1} n_{i} + n_{i-1}$ for $i\ge  2$~\cite[Theorem 1.3]{Bugeaud}, and so $a_{i+1}\le A$ for all $i\ge 2$, which gives (iv). Finally, Bugeaud~\cite[Theorem 1.9 and Definition 1.3]{Bugeaud} gives that (iv) implies (ii), and it is immediate that (ii) implies (i). This  completes the proof.
\end{proof}

We follow common convention, saying that $t$ is \emph{badly approximable} if it satisfies~(i)-(iv) in Proposition~\ref{prop:bad}.  Because of (iv), which we do not directly use here, badly approximable $t$ are sometimes called irrational numbers of \emph{bounded type}.  

Our proof that $\alpha$ is transcendental is substantially simpler if $\theta$ is well (\ie not badly) approximable.  Since the set of all badly approximable numbers is small, having \eg zero Lebesgue measure in $\R$~\cite[Theorem 196]{HW}, it is reasonable to pose the following. 
\begin{question}
Does there exist a Gaussian integer $\zeta$ with argument $2\pi\theta$ for $\theta\in\R$ irrational and well approximable?
\end{question}

\noindent
Unfortunately, the answer is not (as far as we are aware) presently known.  So our arguments will deal with the possibility that $\theta$ is badly approximable, too.
%
%
\subsection{A theorem of Evertse}\label{sec:Evertse}

We now introduce one of our two main technical tools for estimating $\Psi_n(\alpha)$.  Let $K$ be a number field of degree $d:=[K:\mathbb{Q}]$.  Let $M(K)$ denote the set of places of $K$.  Recall (from  e.g.~\cite{EG}) that $M(K)$ is the disjoint union of the set $M_{\rm inf}(K)$ of infinite places and the set $M_{\rm fin}(K)$ of finite places of $K$.  A place $v\in M(K)$ determines a normalized absolute value $|\cdot|_v:K \to [0,\infty)$ as follows.  

If $v\in M_{\rm fin}(K)$ is finite, corresponding to a prime ideal $\mathfrak{p}$ of the ring of integers $\cO_K$ of $K$, then the order ${\rm ord}_\mathfrak{p} x$ of $x\in \cO_K$ is the largest power $m\geq 0$ such that $x \in \mathfrak{p}^m$.  For general $x \in K^\times$, one sets ${\rm ord}_\mathfrak{p} x :=  {\rm ord}_\mathfrak{p} a - {\rm ord}_\mathfrak{p}b$, where $a,b\in\cO_K$ satisfy $x = a/b$.  Then 
$$
|x|_v:= \left\{\begin{array}{ll} 
    0, & \text{if } x=0 \\
    \mathrm{N}(\mathfrak{p})^{-{\rm ord}_\mathfrak{p}(x)}, & \text{if } x\ne 0.    
               \end{array}
\right.
$$ 
where $\mathrm{N}(\mathfrak{p})$ is the cardinality of the finite field $\mathcal{O}_K/\mathfrak{p}$. 
If $v\in M_{\rm inf}(K)$ is an infinite place, then $v$ is either \emph{real} or \emph{complex}. In the first case, $v$ corresponds to a real embedding $\tau\colon K\to \mathbb{R}$, and we take $|x|_v=|\tau(x)|$, where $|\cdot|$ is the ordinary absolute value on $\mathbb{R}$.  In the second case, $v$ corresponds to a conjugate pair $\tau,\bar{\tau}\colon K\to \bC$ of complex embeddings, and we take $|x|_v = |\tau(x)|^{2} = |\bar{\tau}(x)|^{2}$.

A nonzero element $x\in K$ has the property that $|x|_v=1$ for all but finitely many places. With the above normalizations, the \emph{product formula} holds:
\begin{equation} 
\label{eqn:product}
  \prod_{v\in M(K)} |x|_v = 1
  \quad\text{for $x\in K^\times$}.
\end{equation}

If $S\subset M(K)$ is a finite set of places containing all infinite places, then we call 
$\mathcal{O}_{K,S}:=\{a\in K\colon |a|_v\le 1~{\rm for~all~}v\in M(K)\setminus S\}$ the ring of \emph{$S$-integers} in $K$.  Note that if $S = M_{\rm inf}(K)$, then $\cO_{K,S} = \cO_K$ is just the usual ring of integers.
Given a vector $\bfx=(x_1,\ldots ,x_m)\in \mathcal{O}_{K,S}^m$ we set
\begin{equation*}
  H_S(\bfx) = \prod_{v\in S} \max\{|x_1|_v,\ldots, |x_m|_v\}.
\end{equation*}
The following general result of Evertse~\cite{Eve84} (also see~\cite[Proposition 6.2.1]{EG}) on unit equations plays a central role in the sequel.
\begin{theorem}\label{thm:Evertse}
  Let $S\subset M(K)$ be a finite set of places of $K$ containing all infinite places, $m\ge 2$ an integer, and $\epsilon>0$.  There is a constant $c=c(K,S,m,\epsilon)>0$ such that if
  $\bfx=(x_1,\ldots ,x_m)\in \mathcal{O}_{K,S}^m$ and
  $\sum_{k\in I}x_k\ne0$ for every nonempty subset $I\subset\{1,2,\dots,m\}$,
  then for any $v_0\in S$
  \begin{equation*}
    |x_1+\cdots + x_m|_{v_0}
    \ge c\frac{\max\{|x_1|_{v_0},\ldots ,|x_m|_{v_0}\}}
    {{H_S(\bfx)^{\epsilon}\prod_{v\in S} \prod_{k=1}^m |x_k|_v}}.
  \end{equation*}
\end{theorem}

We refer to any quantity of the form $\sum_{k\in I} x_k$, with $I\subset \{1,\dots, m\}$ non-empty, as a \emph{non-trivial subsum} of $x_1+\dots+x_m$.  The assumption that no non-trivial subsum vanishes implies among other things that $x_k\ne 0$ for all $k$.
%
%
\subsection{Initial choices and estimates}\label{sec:estimation}
From now on \emph{we assume that $\alpha$ is an algebraic number}, our final goal being to reach a contradiction.  We fix the number field in the previous subsection to be a(n embedded) Galois extension $K\subset \C$ of $\Q$ that contains $\alpha$, $\bar\alpha$ and $\imunit$.  Any other embedding $K\hookrightarrow\C$ restricts to either the identity or $z\mapsto \bar z$ on $\Q(\imunit)$.  
Hence every infinite place $v$ of $K$ is complex, and the restriction of $|\cdot|_v$ to $\Q(\imunit)$ is the same for all infinite places $v\in M_{\inf}(K)$.  We take $v_0 \in M_{\rm inf}(K)$  to be the infinite place corresponding to the given embedding; \ie $|a|_{v_0} := |a|^2$, where $|\cdot|$ is the restriction to $K$ of the usual absolute value on $\C$.  

We let $\Gamma = \pm\Gamma_0 \cup (\Gamma_0 - \Gamma_0)$.  Then $\Gamma$ contains all coefficients $\gamma(j)$ in the series defining $\Phi(z)$ as well as all differences $\gamma(j)-\gamma(j-n)$, $j>n$.  Specifically, $\Gamma$ is the set of 25 Gaussian integers 
\begin{equation*}
  \Gamma:=\{
  0,\pm2,\pm2\imunit,\pm1\pm2\imunit,
  \pm4,\pm4\imunit,\pm2\pm4\imunit,
  \pm3\pm2\imunit,\pm1\pm4\imunit\};
\end{equation*}
Note for later estimates that if $\gamma\in\Gamma\setminus\{0\}$ and $v\in M_{\inf}(K)$, then
\begin{equation}\label{e303}
  4\leq |\gamma|_v = |\gamma|_{v_0} = |\gamma|^2 \leq 20.
\end{equation}
Finally, we fix
$S\subset M(K)$ to be the set of all infinite places of $K$ together with all finite places $v$ such that $|x|_v \neq 1$ for some $x\in\{\alpha,\bar\alpha\}\cup \Gamma$.  

\begin{lemma}\label{lem:R}
  There is a positive constant $R$ (depending on $\Gamma$ and $\alpha$) such that for any positive integer $n$ and any degree $n$ polynomial 
  $P(z,w) = \sum_{0\leq i+j\leq n} \gamma_{ij}z^i w^j$ with coefficients $\gamma_{ij}\in \Gamma$, the quantity $x = P(\alpha,\bar\alpha)$ satisfies
  \begin{equation*}
      \prod_{v\in S} |x|_v \leq \prod_{v\in M(K)} \max\{|x|_v,1\} \leq R^n.
  \end{equation*}
\end{lemma}

The number $x$ in this lemma is an $S$-integer by construction.  Though the polynomial $P$ used to define $x$ need not be unique, we will be somewhat imprecise and say that \emph{$x$ is a polynomial of degree $n$ in $\alpha$ and $\bar\alpha$ with coefficients in $\Gamma$}.  Whenever we apply Theorem~\ref{thm:Evertse}, it will be to a vector $(x_1,\dots,x_m)$ whose components are all polynomials of this sort.

\begin{proof}
  Pick a positive integer $b$ such that $b\Gamma$, $b\alpha$, and $b\bar{\alpha}$ are all contained in $\mathcal{O}_K$.  Then 
  \begin{equation*}
    b^{n+1} x
    = \sum_{0\le i+j\le n} (b\gamma_{ij})(b\alpha )^i(b\bar\alpha)^j b^{n-i-j} \in \mathcal{O}_K.
  \end{equation*}
  Thus $|b^{n+1}|_v\le 1$ and $|b^{n+1} x|_v\le 1$ for every place $v\in M_{\rm fin}(K)$, so  
\begin{equation*}
    \prod_{v\in M_{\rm fin}(K)} \max\{1,|x|_v\}
    \le 
    \prod_{v\in M_{\rm fin}(K)} |b^{-(n+1)}|_v 
    = \prod_{v\in M_{\rm inf}(K)} |b^{(n+1)}|_v
    = (b^{n+1})^d,
\end{equation*}
  where we used the product formula~\eqref{eqn:product} and the fact that the degree $d$ of $K$ is twice the number of (complex) infinite places.  Let $R_0$ be the maximum of $1$ and the quantities $|\sigma(\alpha)|$ as $\sigma$ ranges over elements of the Galois group $\operatorname{Gal}(K:\bQ)$. 
  Then for any $v\in M_{\inf}(K)$, we have $|x|_v \leq (n+1)^4 R_0^{2n}(\max_{\gamma\in\Gamma}|\gamma|_v) = 20 (n+1)^4R_0^{2n}$, and
\begin{equation*}
  \prod_{v\in M_{\rm inf}(K)} \max\{1,|x|_v\} 
  \leq 20^{d/2} (n+1)^{2d}R_0^{dn}.
\end{equation*}
Putting the estimates for finite and infinite places together then gives
\begin{equation*}
\prod_{v\in S} |x|_v \leq \prod_{v\in M(K)} \max\{1,|x|_v\} \leq (\sqrt{20}(n+1)^2R_0^n)^d b^{(n+1)d} \leq R^n
\end{equation*}
for $R$ (depending on $R_0$, $b$ and $d$) large enough and all $n\geq 1$.
\end{proof}
\begin{corollary}\label{cor:HSestimate}
  If $x_1,\dots,x_m$ are polynomials as in Lemma~\ref{lem:R} and $\sum_{k=1}^m\deg x_k\le n$, then
  \begin{equation*}
    H_S(x_1,\dots,x_m) \leq R^n.
  \end{equation*}
\end{corollary}
\begin{proof}
  Let $n_k = \deg x_k$.  Then by Lemma~\ref{lem:R},
  \begin{equation*}
    H_S(x_1,\dots,x_m)
    = \prod_{v\in S} \max\{|x_1|_v,\dots,|x_m|_v\}
    \leq \prod_{v\in S} \prod_{k=1}^m \max\{1,|x_k|_v\}
    \leq \prod_{k=1}^m R^{n_k}
    \le R^n.
  \end{equation*}
\end{proof}

We conclude by noting that the left-hand estimate in Lemma~\ref{lem:R} can be strengthened when $x$ is a monomial.

\begin{lemma}
\label{lem:monomial}
If $x = \gamma \alpha^j \bar \alpha^k$ for some non-zero $\gamma\in\Gamma$, then $\prod_{v\in S} |x|_v = 1$.
\end{lemma}

\begin{proof}
This follows from the product formula~\eqref{eqn:product} and the fact that $|x|_v = |\gamma|_v|\alpha|_v^j|\bar\alpha|_v^k = 1$ for all places $v\notin S$. 
\end{proof}

%
%
\subsection{The well approximable case}\label{sec:well}
From now on we let $m_i/n_i$, $i\in\N$ denote the continued fraction approximants of $\theta$.  In this section we complete the proof that $\alpha$ is transcendental under the assumption that $\theta$ is well approximable.

\begin{proposition}\label{prop:well}
Suppose that $\theta$ is well approximable.  Then, for any $C\ge 1$, there are arbitrarily large $n\in\N$ such that all indices $j\in (n,Cn]$ are $n$-regular.  
\end{proposition}

\begin{proof}
Let $\epsilon = \frac1{16(C+1)}$.  Since $\theta$ is well-approximable, Proposition~\ref{prop:bad} (i) says that there exist infinitely many $n$ such that $|n\theta-m|<\frac{\epsilon}n$ for some $m\in\N$ coprime to $n$. We claim that any such $n$ will do.
  
To see this, fix $j\in(n,Cn]$ and let $k$ be the integer closest to $8j\theta$.  If $|j\theta - k/8| \geq \frac\epsilon{n}$, then either $j\theta$ and $(j-n)\theta$ are both equivalent mod $1$ to elements of $(\frac{k-1}8,\frac{k}8)$, or both are equivalent to elements of $(\frac{k}8,\frac{k+1}8)$.  Hence (see Figure~\ref{F210}) $\gamma(j\theta) = \gamma((j-n)\theta)$, \ie $j$ is $n$-regular.  If instead $|j\theta - k / 8| < \frac\epsilon{n}$, then
\begin{equation*}
  |8mj - kn| \leq 8|j(m-n\theta) + n(j\theta - k/8)| < 8\epsilon\left(\frac{j}n + 1\right) \leq 8\epsilon (C+1) = \frac12.
\end{equation*}
Hence $8mj = kn$, and since $\gcd(m,n) = 1$, it follows that $j = \frac{k'n}{8}$ where $k' = k/m \in (8,8C]$ is an integer.  Then we have on the one hand that
\begin{equation*}
j\theta - \frac{k}8 = \frac{k'}8(n\theta - m)
\end{equation*}
but on subtracting $n\theta - m$ from both sides, we also obtain
\begin{equation*}
(j-n)\theta - \frac{k-8m}{8} = \frac{k'-8}8(n\theta -m).
\end{equation*}
Since $k'>8$, the right sides of these two equations have the same sign; and their magnitudes are each bounded above by
$
\frac{k'}8|n\theta - m| = |j\theta - \frac{k}8| < \frac1{16}
$
because of our choice of $k$.
So if $n\theta-m > 0$, then $j\theta \in (\frac{k}8,\frac{k+1}8)$ and $(j-n)\theta \in (\frac{k}8- m, \frac{k+1}8 - m)$; and if $n\theta - m< 0$, then $j\theta \in (\frac{k-1}8,\frac{k}8)$ and $(j-n)\theta \in (\frac{k-1}8 - m,\frac{k}8 - m)$.  Either way, $\gamma(j) = \gamma(j-n)$, \ie $j$ is $n$-regular.
\end{proof}

Now let $(x_1,x_2) = (-2|1-\alpha^n|^2,2|1-\alpha^n|^2\re \Phi_n(\alpha))$.

\begin{corollary}
\label{cor:well}
If $\theta$ is well-approximable then for any $C\geq 1$ there are arbitrarily large $n\in\N$ such that
\begin{equation*}
|x_1+x_2| \leq \frac{8\sqrt{5}|\alpha|^{Cn}}{1-|\alpha|}.
\end{equation*}
\end{corollary}

\begin{proof}
From Equation~\eqref{eqn:error} and $\re \Phi(\alpha) = 1$, one sees that
\begin{equation*}
|x_1+x_2| = |\Psi_n(\alpha)| \leq 4\left|\sum_{j>n} (\gamma(j)-\gamma(j-n))\alpha^j\right| = 4\left|\sum_{j>Cn} (\gamma(j)-\gamma(j-n))\alpha^j\right|
\end{equation*}
for any $n$ large enough that Proposition~\ref{prop:well} holds.  The estimate in the corollary now follows from the fact that no element of $\Gamma$ has magnitude larger than $\sqrt{20}$.
\end{proof}

We can apply Theorem~\ref{thm:Evertse} to get a complementary bound for $|x_1+x_2|$.  
Take $v_0$ and $S$ as in the beginning of~\S\ref{sec:estimation}.  Note that
\begin{equation}
\label{eqn:x1x2}
  x_1 = -2(1-\alpha^n)(1-\bar\alpha^n),
  \quad
  x_2 = (1-\bar\alpha^n)\sum_{j=1}^n \gamma(j)\alpha^j
  + (1-\alpha^n)\sum_{j=1}^n\overline{\gamma(j)}\bar\alpha^j
\end{equation}
are polynomials in $\alpha$ and $\bar\alpha$ with degree $2n$ and coefficients in $\Gamma$.
Further, non-trivial subsums of $x_1+x_2$ do not vanish: $|x_1|>1$ for large $n$ because $|\alpha|<1$,
and Lemma~\ref{L118} tells us that $x_2\neq 0$ and $x_1+x_2\ne 0$ for large $n$.
Hence Theorem~\ref{thm:Evertse}, together with Lemma~\ref{lem:R} and Corollary~\ref{cor:HSestimate}, says for any $\epsilon>0$ that
\begin{equation*}
  |x_1 + x_2|^2 = |x_1+x_2|_{v_0}
  \geq c\frac{\max\{|x_1|^2,|x_2|^2\}}{H_S(x_1,x_2)^\epsilon 
    \prod_{v\in S}|x_1|_v|x_2|_v} 
  \geq \frac{c}{R^{4\epsilon n} \cdot R^{4n}}
  = cR^{-4n(1+\epsilon)}.
\end{equation*}
If $\theta$ is well-approximable, then we can compare this lower bound for $|x_1+x_2|$ with the upper bound from Corollary~\ref{cor:well}, obtaining that
\begin{equation*}
|\alpha|^{Cn} \geq c'R^{-2n(1+\epsilon)}
\end{equation*}
for any $C\geq 1$, $\epsilon>0$ fixed, some constant $c'=c'(\epsilon) >0$ and arbitrarily large $n\in\N$. 
Taking $\epsilon = 1$ and $C$ large enough, \eg $C = - 5\frac{\log R}{\log|\alpha|}$, we arrive at a contradiction.  So if $\theta$ is well-approximable then $\alpha$ is transcendental. 
%
%
\subsection{Unit equations}
We need a little extra machinery from the theory of unit equations to deal with the possibility that $\theta$ is badly approximable. Specifically, we need the following result due to Evertse, Schlickewei and Schmidt, see~\cite[Theorem~1.1]{ESS02} and also~\cite[Theorem~6.1.3]{EG}.  To state the theorem, we recall that if $a_1,\dots,a_m\in L$ are (non-zero)
elements of a field $L$, then a solution $y_1,\dots,y_m\in L$ of
\begin{equation*}
a_1y_1+\dots + a_m y_m = 1
\end{equation*}
is called \emph{non-degenerate} if non-trivial subsums of the left side do not vanish.  And a multiplicative subgroup $H\subset (L^\times)^m$ is said to have \emph{rank} $r<\infty$ if there is a free abelian subgroup $H'$ of rank $r$ such that $H/H'$ is finite.

\begin{theorem}\label{thm:unit}
  Let $L$ be a field of characteristic zero,
  $a_1,\dots,a_m\in L^\times$, and $H\subset(L^\times)^m$ a subgroup of finite rank.
  Then there are only finitely many non-degenerate solutions $(y_1,\dots,y_m)\in H$ of the equation
  $a_1y_1+\dots+a_my_m=1$.
\end{theorem}
Note that while Theorem~6.1.3 in~\cite{EG} is only stated for $m\ge 2$, it is also valid (trivially) when $m=1$.
To apply the theorem, let $K$ and $\Gamma$ be as in~\S\ref{sec:estimation}. 

\begin{lemma}\label{L122}
  The numbers $\alpha$ and $\bar\alpha$ generate a
  free multiplicative subgroup of $\bC^\times$.
\end{lemma}
\begin{proof}
  We have $\arg(\alpha)=2\pi\theta$, so if $\alpha^i\bar\alpha^j=1$, then
  $i=j$, as $\theta$ is irrational. But then $\alpha^i\bar\alpha^j=|\alpha|^{2i}$,
  and hence $i=j=0$ since $|\alpha|<1$.
\end{proof}

\begin{corollary}\label{cor:nosubsum}
  For any integer $m\ge1$ there exists $N = N(m)\in\N$ such that 
  \begin{equation*}
    \gamma_1 \alpha^{i_1}\bar\alpha^{j_1}
    + \dots + \gamma_m \alpha^{i_m}\bar\alpha^{j_m} \neq 0
  \end{equation*}
  whenever $\gamma_1,\dots,\gamma_m \in \Gamma$ are not all zero, and $|i_k - i_\ell| + |j_k - j_\ell| \geq N$ for all $k \neq \ell$.
\end{corollary}

\begin{proof}
  We may assume $m\ge2$. Since $\Gamma$ is a finite set, it suffices to consider a
  fixed vector $(\gamma_1,\dots,\gamma_m)$, and we may further assume (after shrinking $m$, if necessary) that $\gamma_k\ne0$ for all $k$.  By Lemma~\ref{L122}, it therefore suffices to prove that for any $(\gamma_1,\dots,\gamma_m)\in(\Gamma\setminus\{0\})^m$
  there are only finitely many non-degenerate solutions $(\alpha^{i'_1}\bar\alpha^{j'_1},\dots,\alpha^{i'_{m-1}}\bar\alpha^{j'_{m-1}})$ to the equation
  \begin{equation*}
    \gamma_1 \alpha^{i'_1}\bar\alpha^{j'_1}
    + \dots
    + \gamma_{m-1} \alpha^{i'_{m-1}}\bar\alpha^{j'_{m-1}}
    + \gamma_m
    =0.
  \end{equation*}
  This follows from Theorem~\ref{thm:unit} with $L=\C$, $m-1$ in place of $m$,
  $a_k=-\gamma_k/\gamma_m$, and $H=G^{m-1}$ where $G\subset\C^\times$ is the multiplicative group generated by $\alpha$ and $\bar\alpha$.
\end{proof}

Theorem~\ref{thm:Evertse} now allows us to render Corollary~\ref{cor:nosubsum} effective:
\begin{corollary}\label{cor:nosmallsubsum}
  Given $\delta,\rho>0$, $C\ge 1$ and an integer $m\ge1$, the following is true for $n$ large enough.  Suppose $i_1,j_1,\dots,i_m,j_m\geq 0$ are integers satisfying 
  \begin{itemize}
  \item $i_k+ j_k \leq Cn$ for all $k$;
  \item $|i_k - i_\ell| + |j_k - j_\ell| \geq \delta n$ for all $k \neq \ell$;
  \end{itemize}
  and suppose $\gamma_1,\dots,\gamma_m\in\Gamma$ do not all vanish.  Then
  \begin{equation}\label{e304}
    |\gamma_1 \alpha^{i_1}\bar\alpha^{j_1}
    + \dots
    +\gamma_m \alpha^{i_m}\bar\alpha^{j_m}|
    \ge |\alpha|^{\min\{i_k+j_k \mid \gamma_k\ne0\}+\rho n}
    \ge |\alpha|^{(C+\rho)n}.
  \end{equation}
\end{corollary}
\begin{proof}
  Suppose without loss of generality that \emph{no} $\gamma_k$ vanishes. Corollary~\ref{cor:nosubsum} tells us that no non-trivial subsum of the sum on the left vanishes.  Let $S\subset M(K)$ and $v_0\in S$ be as in the beginning of~\S\ref{sec:estimation} and $(x_1,\dots,x_m)$ be the vector of monomials 
  $x_k = \gamma_k \alpha^{i_k}\bar\alpha^{j_k}$.  Lemma~\ref{lem:monomial} tells us that $\prod_{v\in S}\prod_{k=1}^m |x_k|_v = 1$.
Further, 
\begin{equation*}
  \max\{|x_1|^2,\ldots ,|x_m|^2\}
  =\max_k|\gamma_k\alpha^{i_k}\bar\alpha^{j_k}|^2
  \geq 4|\alpha|^{2\min_k(i_k+j_k)},
\end{equation*}
using $|\alpha|<1$ and Equation~\eqref{e303}, and Corollary~\ref{cor:HSestimate} gives
\begin{equation*}
  H_S(x_1,\dots,x_m)
  \leq R^{Cmn},
\end{equation*}
for $R>0$ as in Lemma~\ref{lem:R}. Theorem~\ref{thm:Evertse} therefore yields
\begin{equation*}
  |x_1+\dots+x_m|^2 = |x_1+\dots+x_m|_{v_0}
  \geq  \frac{4c|\alpha|^{2\min_k (i_k+j_k)}}{R^{Cmn\epsilon}}.
\end{equation*}
Choosing $\epsilon>0$ small enough that $R^{-Cm\epsilon} > |\alpha|^{2\rho}$ guarantees that the first inequality of Equation~\eqref{e304} holds for large $n$, completing the proof.
\end{proof}
%
%
\subsection{The badly approximable case}\label{sec:badly}
It remains to treat the case when $\theta$ is badly approximable.  We recall (see~\cite[Theorems~167,~171]{HW}) that the continued fraction approximants of an irrational number alternate between over- and under-approximating, \ie if the approximants $m_i/n_i$ of $\theta$ are indexed so that $n_0 = 1$, then we have for any odd index $i$ that $m_{i-1}/n_{i-1} < \theta < m_i/n_i$.  

The next result serves as an alternative to Proposition~\ref{prop:well}.  
\begin{proposition}
\label{prop:badly}
Suppose $\theta$ is badly approximable.  Then there exist $B>0$, $\delta>0$ and arbitrarily large $n\in\N$ such that 
\begin{itemize}
\item[(i)] $j-n\geq \delta n$ for any $n$-irregular index $j>n$;
\item[(ii)] $|j-j'| \geq \delta n$ for any distinct $n$-irregular indices $j,j' > n$;  
\item[(iii)] $|j-j'-n|\geq \delta n$ for any $n$-irregular indices $j,j'> n$ such that $j\neq j'+n$;
\item[(iv)] for any $C\geq 1$, there are at most $C/\delta$ $n$-irregular indices in the interval $(n,Cn]$, and at least one $n$-irregular index in the interval $(Cn,BCn]$.
\end{itemize}
\end{proposition}

\begin{proof}
By hypothesis (see Proposition~\ref{prop:bad}) there exists $\kappa >0$ such that $|n\theta - m| \geq \kappa/n$ for any integers $m,n$ with $n>0$.  In what follows we take $m=m_i$, $n=n_i$ with $i$ odd.

Suppose that $j>n$ is $n$-irregular and let $k$ be the integer closest to $8j\theta$.  Since $m/n$ is a continued fraction approximant of $\theta$, we have $|n\theta - m|<1/n$.  So one can argue as in the second paragraph of the proof of Proposition~\ref{prop:well} to show that $|j\theta - k/8|<\frac1n$.  Hence
\begin{equation*}
  \frac{\kappa}{8(j - n)} \leq |8(j-n)\theta - (k-8m)| \leq 8|j\theta - k/8| + 8|n\theta-m| < \frac{16}n.
\end{equation*}
So $j-n>\frac{\kappa n}{128}$.  And if $j'>n$ is another $n$-irregular index, then $|j'\theta - k'/8|<\frac1n$ for some $k'\in\N$.  Hence if $j'\neq j$,
\begin{equation*}
  \frac{\kappa}{8|j'-j|} \leq |8(j'-j)\theta - (k'-k)| \leq 8|j\theta - k/8| + 8|j'\theta - k'/8| < \frac {16}{n},
\end{equation*}
so $|j'-j| >\frac{\kappa n}{128}$.  Similarly if $j'+n\neq j$, then $|j\theta - k/8|$, $|j'\theta -k'/8|$ and $|n\theta-m|$ are all less than $1/n$, so now the triangle inequality gives
$$
\frac\kappa{8|j-j'-n|} \leq |8(j-j'-n) \theta - (k-k'+8m)| < 24/n,
$$
i.e. $|j-j'-n|>\frac{\kappa n}{192}$. All told,~(i)--(iii) hold with $\delta = \frac{\kappa}{192}$.

The first part of~(iv) follows immediately from~(ii). To prove the second part, pick $m'/n' = m_{i'}/n_{i'}$ to be the continued fraction approximant of $\theta$ with minimal even index $i'$ such that  $n_{i'} > Cn$.  Since $\theta$ is badly approximable, we have $n' \leq A^2Cn$, where $A$ is the constant in the third condition of Proposition~\ref{prop:bad}.  Since $i$ is odd and $i'>i$ is even, we also have
$m'/n' < \theta < m/n$.  Thus
\begin{equation*}
0 < n'\theta-m' < m-n\theta < \frac1n,
\end{equation*}
where the middle inequality comes from the fact that continued fraction approximants of $\theta$ improve as the denominators $n<n'$ increase.
Assuming $n\geq 8$, we infer that $n'\theta$ is equivalent $\mymod 1$ to an element of $(0,\frac18)$.
The inequalities above give
\begin{equation*}
-\frac1n < (n' + n)\theta - (m'+m) < 0
\end{equation*}
so that $(n'+n)\theta$ is equivalent $\mymod 1$ to an element of $(\frac78,1)$.
Then $\gamma(n')=1-2\imunit$ and $\gamma(n+n')=1+2\imunit$, see Figure~\ref{F210},
so the index $j = n' + n$ is $n$-irregular.  Since $Cn < j < (A^2C+1)n$, we may take $B=A^2+1$ to conclude the proof.
\end{proof}

We define $\beta_{i,j}(n)\in\Gamma$ for $i,j\in\N$ by (see Equation~\eqref{eqn:error})
\begin{equation}
\label{e123}
\Psi_n(z) = 2\re\left((1-\bar z^n)\sum_{j>n}(\gamma(j)-\gamma(j-n)) z^j\right)
  = \sum_{i+j>n} \beta_{ij}(n) z^i\bar z^j,
\end{equation}
noting that $\beta_{ij}(n) = \overline{\beta_{ji}(n)}$ is non-zero if and only if one of the indices $i$ or $j$ is
$n$-irregular (hence $>n$) and the other is equal to $0$ or $n$.  Proposition~\ref{prop:badly} implies that for suitable $n$, the indices of non-vanishing $\beta_{ij}(n)$ are well-separated:

\begin{corollary}\label{cor:badly}
  Suppose $\theta$ is badly approximable and let $\delta>0$ be as in Proposition~\ref{prop:badly}. Then, for every $C>\delta/4$ there exists an integer $r\in[0,4C/\delta)$ such that the following assertions hold for infinitely many $n$:
  \begin{itemize}
  \item[(i)]
    if $i,i',j,j'\in\N$ are such that 
    $\beta_{ij}(n)\neq 0$ and $\beta_{i'j'}(n) \neq 0$, then 
    \begin{itemize}
    \item[(a)] $(i,j) = (i',j')$ or $|i-i'| + |j-j'| \geq \delta n$; and
    \item[(b)] $i+j = i'+j'$ or $|(i+j) - (i'+j')| \geq \delta n$;
    \end{itemize}
  \item[(ii)]
    precisely $r$ of the coefficients $\beta_{ij}(n)$ with $i+j\in(n,Cn]$ are non-vanishing.
  \end{itemize}
\end{corollary}

\begin{proof}
  Let $B>0$ be as in Proposition~\ref{prop:badly}, and let $Z\subset\N$ be an infinite subset such that all the assertions of that proposition hold for all $n\in Z$.
  
  Then~(a) follows from Proposition~\ref{prop:badly}~(i)-(ii).  Similarly,~(b) follows from Proposition~\ref{prop:badly}~(i)-(iii): indeed, we can assume, without loss of generality, that $j$ and $j'$ are irregular, and in this case $i-i'\in\{0,n,-n\}$.  

  To prove~(b), set 
  \begin{equation*}
    r_n = \#\{(i,j)\mid \beta_{i,j}(n)\neq 0, i+j\in (n,Cn]\}
  \end{equation*} 
  for any $n\in Z$. 
  In each pair $(i,j)$ being counted, one component is $n$-irregular and the other is equal to either $0$ or $n$.  So Proposition~\ref{prop:badly}~(iv) implies that $r_n <4C/\delta$ for all $n\in Z$.  Hence we can take $r = \liminf r_n$ to be the smallest value of $r_n$ that occurs for infinitely many $n$.  
\end{proof}

Continuing to suppose that $\theta$ is badly approximable, we let $\delta>0$ be as in Proposition~\ref{prop:badly}, and fix $C>\max\{1,\delta/4\}$ (to be specified more precisely below). Let $r\ge 0$ and $n$ be as in Corollary~\ref{cor:badly}. Pick $\rho\in(0,\delta)$.
We will apply Theorem~\ref{thm:Evertse} to the vector
\begin{equation*}
  \bfx = (x_1,x_2,\dots,x_{r+2}) \in \mathcal{O}_{K,S}^{r+2},
\end{equation*}
where 
$$
x_1 = -2|1-\alpha^n|^2, \quad x_2 = 2|1-\alpha^n|^2\re \Phi_n(\alpha),
$$ 
and $x_3,\dots,x_{r+2}$ are the non-vanishing terms $\beta_{ij}(n)\alpha^i\bar\alpha^j$ with $i+j\leq Cn$ in the formula~\eqref{e123} for $\Psi_n(\alpha)$.  From Equation~\eqref{eqn:error} and $\re \Phi(\alpha)=1$ we get
\begin{equation*}
x_1+x_2 = 2|1-\alpha^n|^2(\re \Phi_n(\alpha)-1) = 2|1-\alpha^n|^2\re(\Phi_n(\alpha)-\Phi(\alpha)) = -\Psi_n(\alpha).
\end{equation*}
Together with Equation~\eqref{e123}, this gives
\begin{equation*}
  x_1+\dots+x_{r+2} = -\sum_{i+j>Cn} \beta_{ij}(n) \alpha^i\bar \alpha^j.
\end{equation*} 
Let $p(n)$ denote the maximum value of $i+j$ such that $\beta_{ij}(n)\ne0$ and $i+j\le Cn$. Let $q(n)$ denote the minimum value of $i+j$ such that $\beta_{ij}(n)\ne 0$ and $i+j>Cn$. By Corollary~\ref{cor:badly}~(b) we have $q(n)\ge p(n)+\delta n$.  

Recall that if $\beta_{ij}(n)\ne0$, then the smaller of the indices $i$ and $j$ must equal either $0$ or $n$.  Therefore, for fixed $n$ and $\ell$, there are at most four non-zero $\beta_{ij}(n)$ with $i+j=\ell$. So from the previous equality we estimate 
\begin{equation}\label{e124}
  \left|\sum_{k=1}^{r+2}x_k\right| \leq \sum_{i+j> Cn} |\beta_{ij}(n)| |\alpha|^{i+j}
  \leq 4 \sqrt{20} \sum_{\ell \geq q(n)} |\alpha|^\ell = \frac{8\sqrt{5}|\alpha|^{q(n)}}{1-|\alpha|},
\end{equation}
where the second inequality uses that $\beta_{ij}(n)\in\Gamma$, hence $|\beta_{ij}(n)| \leq \sqrt{20}$.

\begin{lemma}\label{lem:nozerosubsum}
  If $I\subset\{1,2,\dots,r+2\}$ is nonempty, then $\sum_{k\in I}x_k\ne 0$.
\end{lemma}
\begin{proof}
We argue by contradiction, so suppose $\sum_{k\in I}x_k=0$.  By Corollaries~\ref{cor:nosubsum} and~\ref{cor:badly} we cannot have $I\subset\{3,\dots,r+2\}$.  On the other hand, for large $n$, $|x_1|>1$ because $|\alpha|<1$, and $|x_2|>1$ because (additionally) $\re \Phi_n(\alpha) \to 1$.  Finally Lemma~\ref{L118} tells us that $x_1+x_2 \neq 0$ when $n$ is large.  So we cannot have $I\subset\{1,2\}$ either.
  
Since $|x_3+\dots +x_{r+2}|\leq \sqrt{20}r|\alpha|^n < 1<|x_1|,|x_2|$, both $1$ and $2$ belong to $I$ when $n$ is large.  If the complement $J = \{1,\dots,r+2\}\setminus I \subset\{3,\dots,r+2\}$ is non-empty, then Corollaries~\ref{cor:nosmallsubsum} and~\ref{cor:badly} imply that 
\begin{equation*}
\left|\sum_{k=1}^{r+2} x_k\right| = \left|\sum_{k\in J}x_k\right|\ge|\alpha|^{p(n)+\rho n},
\end{equation*}
which contradicts Equation~\eqref{e124} for large $n$ since $p(n) \leq q(n) - \delta n$ and $\rho<\delta$.  

Thus $I=\{1,\dots,r+2\}$, which gives
\begin{equation*}
    0 = -\sum_{j=1}^{r+2} x_j
    =\sum_{i+j > Cn}\beta_{ij}(n)\alpha^i\bar\alpha^j
    =\sum_{i+j=q(n)}\beta_{ij}(n)\alpha^i\bar\alpha^j
    +\sum_{i+j\ge q(n)+\delta n}\beta_{ij}(n)\alpha^i\bar\alpha^j,
\end{equation*}
where we have used Corollary~\ref{cor:badly}~(b) in the last equality.  Hence the two sums on the right have the same magnitude.
Further, $q(n)\le BCn$ by Proposition~\ref{prop:badly}~(iv), so Corollary~\ref{cor:nosmallsubsum} implies
\begin{equation*}
    |\alpha|^{q(n)+\rho n} \leq \left|\sum_{i+j=q(n)}\beta_{ij}(n)\alpha^i\bar\alpha^j\right|
    =
    \left|\sum_{i+j\ge q(n)+\delta n}\beta_{ij}(n)\alpha^i\bar\alpha^j\right|
    \leq 
    \frac{8\sqrt{5}|\alpha|^{q(n)+\delta n}}{1-|\alpha|},
\end{equation*}
where the first inequality follows from Corollary~\ref{cor:nosmallsubsum} and Corollary~\ref{cor:badly}(i), and the second inequality is obtained in the same way as Equation~\eqref{e124}. Since $\rho<\delta$, this is a contradiction for large $n$.
\end{proof}

We are ready to invoke Theorem~\ref{thm:Evertse} one last time, with $S\subset M(K)$ and $v_0\in S$ as in the beginning of~\S\ref{sec:estimation}.  From
Lemma~\ref{lem:monomial} and then Lemma~\ref{lem:R}, we obtain 
\begin{equation*}
  \prod_{v\in S}|x_1|_v\dots |x_{r+2}|_v = \prod_{v\in S} |x_1|_v|x_2|_v \leq R^{4n},
\end{equation*} 
since $x_3,\dots,x_{r+2}$ are monomials in $\alpha,\bar\alpha$, and $x_1,x_2$ are polynomials of degree $2n$, see Equation~\eqref{eqn:x1x2}.  Further, $x_3,\dots,x_{r+2}$ have degree at most $Cn$, so Corollary~\ref{cor:HSestimate} gives
\begin{equation*}
  H_S(\mathbf{x}) \leq R^{(Cr+4)n}.
\end{equation*}

Lemma~\ref{lem:nozerosubsum} says that non-trivial subsums of $x_1+\dots+x_{r+2}$ do not vanish.  So for fixed $\epsilon>0$, Theorem~\ref{thm:Evertse} yields
\begin{equation*}
  |x_1+\dots+x_{r+2}|^2
  \ge c\frac{\max\{|x_1|^2,\dots,|x_{r+2}|^2\}}{H_S(\mathbf{x})^\epsilon 
    \prod_{v\in S}\prod_{k=1}^{r+2} |x_k|_v} 
  \geq \frac{c}{R^{(Cr+4)n\epsilon} \cdot R^{4n}}
  = cR^{-(4 + (Cr+4)\epsilon)n},
\end{equation*}
for large $n$ since $|x_1|\geq 1$.
Using the bound in the other direction from Equation~\eqref{e124}, we infer that if  $n\in\N$ is as in Corollary~\ref{cor:badly} and is large enough, then 
\begin{equation*}
  |\alpha|^{2Cn} \geq |\alpha|^{2q(n)} \geq c'R^{-(4 + (Cr+4)\epsilon)n},
\end{equation*}
for some constant $c' = c'(C,\epsilon)>0$.  So if above we fix $C>1$ such that $|\alpha|^C \le R^{-3}$ and then set $\epsilon = 1/(Cr+4)$, we obtain $R^{-6n}\geq c'R^{-5n}$ for arbitrarily large $n$, which is a
contradiction.  We conclude that if $\theta$ is badly approximable, $\alpha$ is transcendental.

\medskip
This completes the proof of the Main Theorem in the introduction.
%
%
%
%

\end{document}